\documentclass[12pt,reqno]{amsart}

\usepackage{amsfonts,amsmath,amsthm}
\usepackage{amssymb,epsfig}
\usepackage{enumerate} 
\usepackage{subcaption}
\usepackage{nicefrac}
\usepackage{hyperref}

\usepackage{color} 
\definecolor{vert}{rgb}{0,0.6,0}

\usepackage[text={425pt,650pt},centering]{geometry}
\usepackage{graphicx}
\usepackage{epsfig}
\usepackage{tikz,esvect}
\usetikzlibrary{3d,calc,shapes}

\usepackage{caption}
\usepackage{color} 
\definecolor{vert}{rgb}{0,0.6,0}

\numberwithin{figure}{section}

\theoremstyle{plain}
\newtheorem{thm}{Theorem}[section]

\newtheorem{defn}{Definition}

\newtheorem{ex}{Example}
\newtheorem{lem}[thm]{Lemma}
\newtheorem{cor}[thm]{Corollary}
\newtheorem{prop}[thm]{Proposition}
\theoremstyle{remark}
\newtheorem{rem}{\bf{Remark}}
\numberwithin{equation}{section}

\DeclareMathOperator*{\essinf}{ess\,inf}

\newcommand{\N}{\mathbb{N}}

\newcommand{\R}{\mathbb{R}}

\newcommand{\T}{\mathbb{T}}
\newcommand{\Z}{\mathbb{Z}}

\newcommand{\Lip}{{\rm Lip\,}}

\newcommand{\al}{\alpha}
\newcommand{\gam}{\gamma}

\newcommand{\ep}{\varepsilon}

\newcommand{\Del}{\Delta}

\newcommand{\ol}{\overline}

\DeclareMathOperator\erf{erf}

\begin{document}

\title[Optimal homogenization rate of viscous HJ equation]
{Optimal rate of convergence in periodic homogenization of viscous Hamilton-Jacobi equations}

\author[J. QIAN, T. SPREKELER, H. V. TRAN, Y. YU]
{Jianliang Qian, Timo Sprekeler, Hung V. Tran, Yifeng Yu}

\thanks{
The work of JQ is partially supported by NSF grants 2012046, 2152011, and 2309534. 
The work of HT is partially supported by NSF CAREER grant DMS-1843320 and a Vilas Faculty Early-Career Investigator Award.
The work of YY is partially supported by NSF grant 2000191.
}

\address[J. Qian]
{
Department of Mathematics, Michigan State University, East Lansing, MI 48824, USA}
\email{jqian@msu.edu}

\address[T. Sprekeler]
{
Department of Mathematics,
National University of Singapore, 
10 Lower Kent Ridge Road, 
Singapore 119076}
\email{timo.sprekeler@nus.edu.sg}

\address[H. V. Tran]
{
Department of Mathematics, 
University of Wisconsin Madison, Van Vleck Hall, 480 Lincoln Drive, Madison, Wisconsin 53706, USA}
\email{hung@math.wisc.edu}

\address[Y. Yu]
{
Department of Mathematics, 
University of California at Irvine, 
California 92697, USA}
\email{yyu1@math.uci.edu}

%\date{\today}

\maketitle

\begin{abstract}
We study the optimal rate of convergence in periodic homogenization of the viscous Hamilton-Jacobi equation $u^\ep_t + H(\frac{x}{\ep},Du^\ep) = \ep \Delta u^\ep$ in $\R^n\times (0,\infty)$ subject to a given initial datum. We prove that $\|u^\ep-u\|_{L^\infty(\R^n \times [0,T])} \leq C(1+T) \sqrt{\ep}$ for any given $T>0$, where $u$ is the viscosity solution of the effective problem. Moreover, we show that the $O(\sqrt{\ep})$ rate is optimal for a natural class of $H$ and a Lipschitz continuous initial datum, both theoretically and through numerical experiments.  It remains an interesting question to investigate whether the convergence rate can be improved when $H$ is uniformly convex. Finally, we propose a numerical scheme for the approximation of the effective Hamiltonian based on a finite element approximation of approximate corrector problems.
\end{abstract}

\section{Introduction}
\subsection{Settings}
For each $\ep>0$, let $u^\ep \in C(\R^n \times [0,\infty))$ be the viscosity solution to 
\begin{equation}\label{eq:C-ep}
\begin{cases}
 u_t^\ep+H\left(\frac{x}{\ep},Du^\ep\right)=\ep \Delta u^\ep \qquad &\text{in} \ \R^n \times (0,\infty),\\
u^\ep(x,0)=g(x) \qquad &\text{on} \ \R^n.
\end{cases} 
\end{equation}
Here, $g\in C^{0,1}(\R^n)$ is a given initial datum and $H = H(y,p)\in \Lip_{\rm loc}(\R^n\times \R^n)$ is a given Hamiltonian that is $\Z^n$-periodic in its $y$-variable and satisfies 
\begin{align}\label{coe}
\essinf_{y \in \R^n} \left\{|H(y,p)|^2+(n+1)D_yH(y,p)\cdot p \right\} \longrightarrow \infty\quad\text{as }|p|\rightarrow \infty.
\end{align}
Then, it is known that $u^\ep$ converges to $u\in C(\R^n \times [0,\infty))$ locally uniformly on $\R^n \times [0,\infty)$ as $\ep \to 0^+$, where $u$ is the viscosity solution to the effective problem
\begin{equation}\label{eq:C} 
\begin{cases}
 u_t+\ol{H}\left(Du\right)=0 \qquad &\text{in} \ \R^n \times (0,\infty),\\
u(x,0)=g(x) \qquad &\text{on} \ \R^n; 
\end{cases} 
\end{equation}
see \cite{LPV, Ev1}. Here, the effective Hamiltonian $\ol{H} \in C(\R^n)$ is determined by $H$ in a nonlinear way through cell problems.  
It is worth noting that if $H = H(y,p)$ is independent of $y$, that is, $H(y,p)=F(p)$, then \eqref{eq:C-ep} becomes the usual vanishing viscosity problem
\begin{equation}\label{eq:F-ep}
\begin{cases}
 u_t^\ep+F\left(Du^\ep\right)=\ep \Delta u^\ep \qquad &\text{in} \ \R^n \times (0,\infty),\\
u^\ep(x,0)=g(x) \qquad &\text{on} \ \R^n,
\end{cases} 
\end{equation}
in which case we have $\ol H = F$.
Both \eqref{eq:C-ep} and \eqref{eq:F-ep} are basic and fundamentally important problems in the theory of viscosity solutions.

Introducing the notation $\T^n:=\R^n/\Z^n$, we now give a precise definition of $\ol H$.
\begin{defn}[Effective Hamiltonian]\label{def: Hbar}
Assume {\rm (A1)--(A2)}.
For each $p \in \R^n$, there exists a unique constant $\ol H(p)\in \R$ such that the cell (ergodic) problem 
\begin{equation} \label{eq:E-p}
 \qquad H(y,p+Dv)=\ol H(p) +\Delta v \qquad \text{for} \ y\in\T^n
\end{equation}
has a continuous viscosity solution $v$. If needed, we write $v=v(y,p)$ or $v=v_p(y)$ to clearly demonstrate the nonlinear dependence of $v$ on $p$. 
In the literature, $v(\cdot,p)$ is often called a corrector.
It is worth mentioning that $v(\cdot,p)$ is unique up to additive constants.
\end{defn}
From now on, we normalize the corrector $v$ so that $v(0,p)=0$ for all $p\in \R^n$.
In fact, $v(\cdot,p)\in C^2(\T^n)$ and $p \mapsto v(\cdot,p)$ is locally Lipschitz. 
Further, the effective Hamiltonian $\ol H$ is locally Lipschitz.

Our main goal in this paper is to obtain the optimal rate for the convergence of $u^\ep$ to $u$, that is, an optimal bound for $\|u^\ep-u\|_{L^\infty(\R^n \times [0,T])}$ for any given $T>0$ as $\ep \to 0^+$.  
Heuristically, thanks to the two-scale asymptotic expansion, 
\begin{equation}\label{expansion}
u^{\ep}(x,t) \approx u(x,t)+\ep v\left(\frac{x}{\ep}, Du(x,t)\right)+O(\ep^2).
\end{equation}
However, this is just a formal local expansion, and it is not clear at all how to obtain the optimal global bound in the $L^\infty$-norm from this.

\subsection{Main results}
We now describe our main results. Let us introduce the set of assumptions (A1)--(A3) given by
\begin{itemize}
\item[(A1)] $H\in \Lip_{\rm loc}(\R^n\times \R^n)$, and $H(\cdot,p)$ is $\Z^n$-periodic for each $p\in \R^n$;

\item[(A2)] $H$ satisfies \eqref{coe};

\item[(A3)] $g\in \Lip(\R^n)$ with $\|g\|_{C^{0,1}(\R^n)}<\infty$.
\end{itemize}
\begin{thm} \label{thm:nd}
Assume {\rm (A1)--(A3)} and fix $T>0$. Then, there exists a constant $C>0$ depending only on $H$, $n$, and $\|g\|_{C^{0,1}(\R^n)}$ such that for $\ep \in (0,1)$ there holds
\[
\|u^\ep-u\|_{L^\infty(\R^n \times [0,T])} \leq C(1+T) \sqrt{\ep},
\]
where $u^\ep$ and $u$ denote the viscosity solutions to \eqref{eq:C-ep} and \eqref{eq:C}, respectively.  
\end{thm}

The above rate $O(\sqrt{\ep})$ turns out to be optimal in the sense that  there exist particular choices of $H$ and $g$ satisfying (A1)--(A3) such that the convergence rate is exactly $O(\sqrt{\ep})$.  Quantitative homogenization for Hamilton-Jacobi equations in the periodic setting has received quite a lot of attention in the past twenty years. The convergence rate $O(\ep^{\nicefrac{1}{3}})$ was obtained for first-order equations first in \cite{CDI}.
In \cite{CCM}, the authors generalized the method in \cite{CDI} to get the same convergence rate $O(\ep^{\nicefrac{1}{3}})$ for the viscous case considered in this paper.
For weakly coupled systems of first-order equations, see \cite{MT}.
For other related works, see the references in \cite{CDI, CCM, MT}.
Of course, the rate $O(\ep^{\nicefrac{1}{3}})$ is not known to be optimal in general.

The optimal rate of convergence $O(\ep)$ for convex first-order equations was recently obtained in \cite{TY}. Moreover,  we expect  that for any given uniformly convex $H$, the convergence rate is $O(\ep)$ for \eqref{eq:C-ep} for generic initial data, which is stronger than the notion of optimality in this paper. We refer to \cite{HJ} for the multi-scale setting.
For earlier progress in this direction with nearly optimal rates of convergence, we refer the reader to  \cite{MTY, Tu, JTY, Coop} and the references therein.
To date, optimal rates of convergence for general nonconvex first-order cases have not been established.

To the best of our knowledge, the optimal rate of convergence for periodic homogenization of viscous Hamilton-Jacobi equations has not been obtained in the current literature.
The rate $O(\ep^{\nicefrac{1}{3}})$ was obtained in \cite{CDI, CCM} by using the doubling variable technique, the perturbed test function method \cite{Ev1}, and the approximate cell problems.
The usage of the approximate cell problems introduces another parameter in the analysis, and as a result, the rate $O(\ep^{\nicefrac{1}{3}})$ was the best one can obtain through this route by optimizing over all parameters.

In this paper, we are able to obtain the $O(\sqrt{\ep})$ convergence rate by dealing directly with the correctors.
A key point is that after normalizing $v(0,p)=0$, we have that $v(\cdot,p)$ is unique, and $p\mapsto v(\cdot,p)$ is locally Lipschitz.
It is worth noting that we do not require convexity of the Hamiltonian in Theorem \ref{thm:nd}.

Here, we will use $H(y,p)=F(p)$ for some choices of  nonlinear $F$ to construct computable sharp examples.  Similar results were known for linear $F$ in the context of conservation laws  \cite{TT1}. The connection between scalar conservation laws and Hamilton-Jacobi equations  is well known to experts. 
 Precisely speaking, in one dimension, if $u=u(x,t)$ is a viscosity solution to $u_t+F(u_x)=0$, then $v=u_x$ is an entropy solution to $v_t+(F(v))_x=0$.   
 The convergence rate of vanishing viscosity in scalar conservation laws has been well studied and the  convergence rate of $O(\sqrt{\ep})$ was known under suitable assumptions \cite{K1976}. 

\begin{thm}\label{thm:new}
Let $n=1$.
Let $F \in \Lip_{\rm loc}(\R)$ be such that 
\[
\begin{cases}
F(p)=p \qquad &\text{ for $p\in [0,1]$}, \\
F(p) \leq p \qquad &\text{ for $p\in [-1,0]$},
\end{cases}
\]
and suppose that $g(x)=\max\{1-|x|,0\}$ for $x\in \R$. Then, for any $\ep \in (0,\frac{1}{4})$ there holds
\[
|u^\ep(0,1)-u(0,1)|  \geq \frac{\mathrm{e}-1}{\sqrt{\pi} \mathrm{e}} \sqrt{\ep},
\]
where $u^\ep$ denotes the viscosity solution to  \eqref{eq:F-ep} and $u$ denotes the viscosity solution to \eqref{eq:C} with $\ol H = F$.
\end{thm}

 We would like to point out that the above $g$ can be replaced by a smooth function (Remark \ref{rem:C2}). Also,   the proof of Theorem \ref{thm:new} leads to the following corollary.

\begin{cor} \label{cor:gen}
Let $n=1$.
Assume that $F \in \Lip_{\rm loc}(\R)$ and that $F$ is linear in $(a,b) \subset \R$ for some given $a<b$.
Then, there exists an initial datum $g \in \Lip(\R)$ such that for any $\ep \in (0,\frac{1}{4})$ we have that
\[
|u^\ep(0,1)-u(0,1)|  \geq c_0  \sqrt{\ep}
\]
for some constant $c_0>0$ depending only on $F$ and $g$, where $u^\ep$ denotes the viscosity solution to  \eqref{eq:F-ep} and $u$ denotes the viscosity solution to \eqref{eq:C} with $\ol H = F$.

\end{cor}

It is also straightforward to generalize Theorem \ref{thm:new} to any dimension in the corollary below, whose proof is essentially the same as that of Theorem \ref{thm:new}.

\begin{cor}\label{cor:new}
Let $F \in \Lip_{\rm loc}(\R^n)$ be such that 
\[
\begin{cases}
F(s e_1)=s \qquad &\text{ for $s\in [0,1]$}, \\
F(se_1) \leq s \qquad &\text{ for $s\in [-1,0]$},
\end{cases}
\]
and suppose that $g(x)=\max\{1-|x_1|,0\}$ for $x\in \R^n$. Then, for any $\ep \in (0,\frac{1}{4})$ there holds
\[
|u^\ep(0,1)-u(0,1)|  \geq \frac{\mathrm{e}-1}{\sqrt{\pi} \mathrm{e}} \sqrt{\ep},
\]
where $u^\ep$ denotes the viscosity solution to  \eqref{eq:F-ep} and $u$ denotes the viscosity solution to \eqref{eq:C} with $\ol H = F$.
\end{cor}

The bound $O(\sqrt{T\ep})$ for $\|u^\ep-u\|_{L^\infty(\R^n \times [0,T])}$  for the vanishing viscosity process of \eqref{eq:F-ep} was obtained in \cite{Fl, CL, Ev2}.
In this situation, we only need to assume that $F$ is locally Lipschitz on $\R^n$ and $g$ is bounded and Lipschitz on $\R^n$ (see e.g., \cite[Theorem 5.1]{CL}).
For the static cases, see \cite{Tr1, Tran}.

Thus, the results of Theorem  \ref{thm:new} and Corollaries \ref{cor:gen}--\ref{cor:new} confirm both the optimality of the convergence rate of  the vanishing viscosity process of \eqref{eq:F-ep} with optimal conditions, and  the optimality of the $O(\sqrt{\ep})$ bound in Theorem \ref{thm:nd}.
See Remark \ref{rem:C2} for Theorem  \ref{thm:new} with a $C^2$ initial condition for each $\ep>0$.
Besides, we provide a generalization of Theorem  \ref{thm:new} in Proposition \ref{prop:generalization} in which for each fixed $\ep\in (0,\frac{1}{4})$, the Hamiltonian $F$ needs not to be linear in any interval in one dimension at the price of  nonconvexity.  Note also that in Corollary \ref{cor:new}, $F$ does not need to be linear in any open set in multiple dimensions.

 Note that all the Hamiltonians in Theorem  \ref{thm:new} and Corollaries \ref{cor:gen}--\ref{cor:new} are not strictly convex and do not have $y$-dependence (i.e., no homogenization effect  is involved).   Hence, it is natural to ask  (I) whether the convergence rate can be improved for strictly/uniformly convex  $H$ and (II) how the $y$-dependence impacts the convergence rate.

(I) has been investigated in the context of one dimensional conservation laws for the vanishing viscosity process of \eqref{eq:F-ep}.  It was proved that the convergence rate can be improved to $O(\ep|\log\ep|)$ for uniformly convex $F$ under some technical assumptions \cite{TT}.  In Section \ref{subsec:quadHam},  we demonstrate this fact for the quadratic Hamiltonian  $F(p)=\frac{1}{2}|p|^2$ in any dimension for general Lipschitz continuous initial data. More  interestingly,  we showed that  for any $C^2$ initial datum $g$, the convergence rate is $O(\ep)$ for a.e. $(x,t)\in \R^n\times (0,\infty)$.   For strictly but not uniformly convex $F$, numerical computation shows that the convergence rate could be various fractions. For instance,  for $F(p)=\frac{1}{4} |p|^4$ in Example \ref{ex:5}, the rate of convergence  for the vanishing viscosity process of \eqref{eq:F-ep} seems to be  $O(\ep^{2/3})$.  This suggests that there might be  a variety of rates $O(\ep^s)$ for $\frac{1}{2}\leq s \leq 1$ for \eqref{eq:F-ep}, which is a new phenomenon.  It will be an interesting project to find an example where a convergence rate $\alpha\in (\frac{1}{2},1)$ can be  established rigorously.

As for (II), it is quite challenging to conduct a theoretical analysis beyond Theorem  \ref{thm:nd} when $y$ is present. In this paper, we will focus on numerical computations to get some rough ideas and inspire interested readers to work on this  subject. Our numerical Examples \ref{ex:10} and \ref{ex:11} show that when $H = H(y,p)$  is strictly convex in $p$ and smooth in $y$, the convergence rate is similar to $O(\ep)$ or $O(\ep|\log \ep|)$. Meanwhile, when the regularity in $y$ is merely Lipschitz continuity, the convergence rate seems to be reduced; see Examples \ref{ex:6} -\ref{ex:9}.

Finally, we discuss the construction of numerical methods for the approximation of the effective Hamiltonian $\ol H$. In particular, we provide a simple scheme to approximate $\ol H$ at a fixed point based on a finite element approximation of approximate corrector problems. 
For related work on the numerical approximation of effective Hamiltonians we refer to \cite{ACC,FR,GLQ,GO,LYZ,Qian} for first-order Hamilton-Jacobi equations without viscosity term, and to \cite{GSS,KS} for second-order Hamilton-Jacobi-Bellman and Isaacs equations.

\subsection*{Organization of the paper}
In Section \ref{sec:simplification}, we use a priori estimates to simplify the settings of the problems.
The proof of the bound in Theorem \ref{thm:nd} is given in Section \ref{sec:proof-main}.
In Section \ref{sec:optimality}, we consider \eqref{eq:F-ep} with  various choices of $F$ and $g$, and obtain the optimality of the bound in Theorem \ref{thm:nd}.
In particular, this section includes a proof of Theorem \ref{thm:new}.
Numerical results for both \eqref{eq:C-ep} and \eqref{eq:F-ep} are studied in Section \ref{sec:num}. The approximation of the effective Hamiltonian is studied in Section \ref{Sec: Num effHbar}.

\section{Settings and simplifications} \label{sec:simplification}
Assume {\rm (A1)--(A3)}.
For $\ep \in (0,1)$, let $u^\ep$ denote the viscosity solution to \eqref{eq:C-ep}. Let $u$ denote the viscosity solution to \eqref{eq:C}. By the comparison principle, we have that
\begin{equation}\label{u-bdd}
\|u_t\|_{L^{\infty}(\R^n \times [0,\infty))}+\|Du\|_{L^{\infty}(\R^n \times [0,\infty))}\leq M
\end{equation}
for $M:=R_0+\max_{\{|p|\leq R_0\}}|\overline H(p)|$, where $R_0:= \|Dg\|_{L^{\infty}(\R^n)}$.

Let us further  assume that 
\begin{equation}\label{u-ep-bdd}
\|u^{\ep}_t\|_{L^{\infty}(\R^n \times [0,\infty))}+\|Du^{\ep}\|_{L^{\infty}(\R^n \times [0,\infty))}\leq C_0
\end{equation}
for a constant $C_0\geq M$ that is independent of $\ep$. Note that  \eqref{u-ep-bdd}  is satisfied if $g\in C^2(\R^n)$ with $\|g\|_{C^2(\R^n)}<\infty$, using the classical Bernstein method based on (A1)--(A2).  Since $g$ is merely assumed to be in $C^{0,1}(\R^n)$ in Theorem \ref{thm:nd},  we will employ a suitable mollification of $g$ in Section \ref{subsec: Removal}  to remove the assumption  \eqref{u-ep-bdd}. 

Accordingly,  values of $H(y,p)$ for $|p|>C_0$ are irrelevant. 
Indeed, letting $\xi \in C^\infty(\R^n,[0,1])$ be a cut-off function satisfying
\[
\xi(p)=1 \text{ if } |p| \leq C_0+1,\qquad
\xi(p)=0 \text{ if } |p| \geq 2(C_0+1),
\]
and introducing
\[
\widetilde H(y,p) := \xi(p) H(y,p) + (1-\xi(p)) |p|^2 \quad \text{ for } (y,p) \in \T^n \times \R^n,
\]
we have that $\widetilde H$ satisfies (A1)--(A2) and $u^\ep$ solves \eqref{eq:C-ep} with $\widetilde H$ in place of $H$.
Therefore, from now on, we can assume that $H$ takes the form of $\widetilde H$, that is, $H$ satisfies
\begin{itemize}
\item[(A4)] $H(y,p) = |p|^2$ for $y\in \T^n$ and $|p| \geq 2(C_0+1)$.
\end{itemize}
Assumption (A4) helps us simplify the situation quite a bit as follows.
For $|p| \geq 2(C_0+1)$, it is clear that $v(\cdot,p)\equiv 0$ and $\ol H(p) = |p|^2$.
Hence, we obtain that $p \mapsto v(\cdot,p)$ is bounded and globally Lipschitz, that is, there exists $C>0$ such that
\begin{equation}\label{v-p-Lip}
\|v(\cdot,p)\|_{L^\infty(\T^n)} \leq C,\quad \|v(\cdot, p) - v(\cdot, \tilde{p})\|_{L^\infty(\T^n)} \leq C |p-\tilde{p}| \qquad \forall p, \tilde{p} \in \R^n.
\end{equation}

If a function $h:\R^n\to\R$ is $\Z^n$-periodic, we can think of $h$ as a function from $\T^n$ to $\R$ as well, and vise versa.
In this paper, we switch freely between the two interpretations.

\section{Proof of Theorem \ref{thm:nd}}\label{sec:proof-main}

\subsection{Part 1: Proof based on \eqref{u-ep-bdd}}\label{subsec: Pf based on} Assumptions (A1)--(A4) are always in force in this section. Let $T>0$ be fixed.
Our goal is to show that there exists a constant $C>0$ depending only on $\|H\|_{C^{0,1}\left(\T^n\times \overline {B(0,2(C_0+1))}\right)}$, $n$, and $C_0$  from (\ref{u-ep-bdd}) such that for any $\ep\in (0,1)$ there holds
\begin{equation}\label{error-bound}
\|u^\ep-u\|_{L^\infty(\R^n \times [0,T])} \leq C(1+T) \sqrt{\ep}.
\end{equation}
The approach here is inspired by that in \cite[Theorem 4.40]{Tran}. We first show that
\begin{align}\label{upperbd}
u^\ep(x,t) - u(x,t) \leq C (1+T) \sqrt{\ep}\qquad \forall (x,t) \in \R^n \times [0,T].
\end{align}
\begin{proof}[Proof of \eqref{upperbd}]
We divide the proof into several steps.

\medskip
\noindent\textbf{Step 0:} We write $\mathcal{A}:=\R^n\times \R^n\times \R^n\times [0,\infty)\times [0,\infty)$. For $K>0$ to be chosen and $\gamma \in (0,\frac{1}{2})$, we introduce the auxiliary function $\Phi_1:\mathcal{A}\rightarrow \R$ given by
\begin{align*}
\Phi_1(a) := u^\ep(x,t) - u(y,s) - \ep v\left(\frac{x}{\ep},\frac{z-y}{\sqrt{\ep}}\right) - \omega(a)\quad\text{for }a = (x,y,z,t,s)\in \mathcal{A},
\end{align*}
where
\begin{align}\label{omega def}
\omega(x,y,z,t,s):= \frac{|x-y|^2 + |x-z|^2 + |t-s|^2}{2 \sqrt{\ep}} + K(t+s) + \gam \sqrt{1+|x|^2}.
\end{align}
Note that there exist $\hat x, \hat y, \hat z\in \R^n$ and $\hat t, \hat s\in [0,\infty)$ such that $\Phi_1$ has a global maximum at the point $\hat{a}:=(\hat x, \hat y, \hat z, \hat t, \hat s)\in \mathcal{A}$. 
We fix such a choice of $\hat{a}$ and introduce $\Phi:\mathcal{A}\rightarrow \R$ given by 
\begin{align*}
\Phi(a):=\Phi_1(a) - \gamma\frac{|a - \hat{a}|^2}{2}\quad\text{for }a\in \mathcal{A}.
\end{align*}
Observe that $\Phi$ has a strict global maximum at the point $\hat{a}$.

\medskip 
\noindent \textbf{Step 1:} We show that 
\begin{align}\label{Goal Step1}
|\hat x - \hat z| \leq C \ep,\qquad |\hat x - \hat y| + |\hat y - \hat z| \leq C \sqrt{\ep},\qquad |\hat t - \hat s| \leq C(1+K) \sqrt{\ep}.
\end{align}
To this end, we first use that $\Phi(\hat a) \geq \Phi(\hat x, \hat y, \hat x, \hat t, \hat s)$ and \eqref{v-p-Lip} to obtain 
\[
(1-\gamma\sqrt{\ep}) \frac{|\hat x-\hat z|^2}{2\sqrt{\ep}} \leq  \ep \left[v\left(\frac{\hat x}{\ep},\frac{\hat x -\hat y}{\sqrt{\ep}}\right) - v\left(\frac{\hat x}{\ep},\frac{\hat z -\hat y}{\sqrt{\ep}}\right) \right] \leq C \sqrt{\ep} |\hat x - \hat z|,
\]
which yields $ |\hat x - \hat z| \leq C \ep$  as $\gamma\sqrt{\ep}\leq\frac{1}{2}$. Then, we use that $\Phi(\hat a) \geq \Phi(\hat x, \hat x, \hat x, \hat t, \hat s)$, \eqref{u-bdd}, and \eqref{v-p-Lip} to find that
\begin{align*}
(1-\gamma\sqrt{\ep})\frac{|\hat x-\hat y|^2+|\hat x-\hat z|^2}{2\sqrt{\ep}} &\leq u(\hat x, \hat s ) - u(\hat y, \hat s) + \ep \left[v\left(\frac{\hat x}{\ep},0\right) - v\left(\frac{\hat x}{\ep},\frac{\hat z -\hat y}{\sqrt{\ep}}\right) \right] \\
& \leq C|\hat x - \hat y| + C \sqrt{\ep} |\hat y - \hat z| \\
&\leq C|\hat x - \hat y| + C \ep^{\nicefrac{3}{2}},
\end{align*}
which yields $|\hat x - \hat y| \leq C \sqrt{\ep}$. Finally, using $\Phi(\hat a) \geq \Phi(\hat x, \hat y, \hat z, \hat t, \hat t)$  and \eqref{u-bdd}, we find
\begin{align*}
(1-\gamma\sqrt{\ep})\frac{|\hat t-\hat s|^2}{ 2\sqrt{\ep}}\leq u(\hat y,\hat t)-u(\hat y,\hat s) + K(\hat t - \hat s)\leq (C+K)|\hat t -\hat s|,
\end{align*}
which yields $|\hat t - \hat s| \leq C(1+K) \sqrt{\ep}$.

\medskip 
\noindent \textbf{Step 2:} For the case $\hat t,\hat s > 0$, we show that
\begin{equation}\label{Goal Step2}
K+ \frac{\hat t - \hat s}{\sqrt{\ep}} + \ol H\left(\frac{\hat z - \hat y}{\sqrt{\ep}} \right) \leq  C \sqrt{\ep} + C \gam.
\end{equation}
Introducing $\varphi:\R^n\times [0,\infty)\rightarrow \R$ defined by $\varphi(x,t):=u^{\ep}(x,t)-\Phi(x,\hat y,\hat z,t,\hat s)$, we see that $u^{\ep}-\varphi$ has a global maximum at $(\hat x, \hat t)$. We compute $\varphi_t(\hat x,\hat t) = K+\frac{\hat t-\hat s}{\sqrt{\ep}}$ and
\begin{align*}
D\varphi(\hat x,\hat t) &= Dv\left(\frac{\hat x}{\ep},\frac{\hat z - \hat y}{\sqrt{\ep}}\right) + \frac{(\hat x - \hat y)+(\hat x - \hat z)}{\sqrt{\ep}}+\gam \frac{\hat x}{\sqrt{1+|\hat x|^2}},\\
\ep \Delta \varphi(\hat x,\hat t) &= \Delta v\left(\frac{\hat x}{\ep},\frac{\hat z - \hat y}{\sqrt{\ep}}\right) + 2n \sqrt{\ep} + \gam \ep \frac{n + (n-1)|\hat x|^2}{(1+|\hat x|^2)^{\nicefrac{3}{2}}}+\gamma n\ep.
\end{align*}
Writing $y_0:=\frac{\hat x}{\ep}$ and $p_0 := \frac{\hat z - \hat y}{\sqrt{\ep}}$, we can use the viscosity subsolution test, \eqref{eq:E-p}, \eqref{Goal Step1}, and local Lipschitz continuity of $H$ to find that
\begin{align*}
&K + \frac{\hat t-\hat s}{\sqrt{\ep}} + \ol H\left(\frac{\hat z - \hat y}{\sqrt{\ep}} \right) = \varphi_t(\hat x,\hat t) + \ol H(p_0) \\&\leq \left[\ep \Delta \varphi(\hat x,\hat t) -\Delta v(y_0,p_0)\right] +  \left[H(y_0,p_0 + Dv(y_0,p_0 ) )- H(y_0,D\varphi(\hat x,\hat t))\right]\\
&\leq  2n\sqrt{\ep} + \gam \ep \frac{n + (n-1)|\hat x|^2}{(1+|\hat x|^2)^{\nicefrac{3}{2}}} + \gamma n\ep+  C\left\lvert 2\frac{\hat z - \hat x}{\sqrt{\ep}}-\gam \frac{\hat x}{\sqrt{1+|\hat x|^2}}  \right\rvert
\\
&\leq C\sqrt{\ep} + C \gamma; 
\end{align*}
i.e., \eqref{Goal Step2} holds. 

\medskip 
\noindent \textbf{Step 3:} For the case $\hat t,\hat s>0$, we show that
\begin{equation}\label{Goal Step3}
K- \frac{\hat t - \hat s}{\sqrt{\ep}} - \ol H\left(\frac{\hat x - \hat y}{\sqrt{\ep}} \right) \leq  C \sqrt{\ep}.
\end{equation}
For $\al>0$, we introduce the auxiliary function $\Psi:\R^n\times \R^n\times [0,\infty)\rightarrow \R$ given by
\begin{align*}
\Psi(y,\xi,s) &:=  u(y,s) + \ep v \left(\frac{\hat x}{\ep}, \frac{\hat z - \xi}{\sqrt{\ep}} \right) + \frac{|\hat x - y|^2 + |\hat t - s|^2}{2 \sqrt{\ep}} + \frac{|y-\xi|^2}{2\al}+Ks \\ &\quad+\gamma \frac{|\hat{y}-y|^2 + |\hat{s}-s|^2}{2}.
\end{align*}
Note that there exist $y_\al,\xi_\al\in \R^n$ and $s_\al\in [0,\infty)$ such that the function $\Psi$ has a global minimum at the point $(y_\al, \xi_\al,s_\al)$. 

\noindent
\textit{Step 3.1:} We first show that
\begin{align}\label{Step3.1}
|y_\al - \xi_\al| \leq C \al\sqrt{\ep},\quad |\hat x - y_\al|\leq C\sqrt{1+\alpha}\sqrt{\ep},\quad |\hat{t}-s_{\alpha}|\leq C(1+K)\sqrt{\ep}.
\end{align}
We use $\Psi(y_\al, \xi_\al,s_\al) \leq \Psi(y_\al, y_\al,s_\al)$  and \eqref{v-p-Lip} to obtain
\[
\frac{|y_\al-\xi_\al|^2}{2\al} \leq \ep \left[  v \left(\frac{\hat x}{\ep}, \frac{\hat z - y_\al}{\sqrt{\ep}} \right) - v \left(\frac{\hat x}{\ep}, \frac{\hat z - \xi_\al}{\sqrt{\ep}} \right)\right] \leq C \sqrt{\ep} |y_\al - \xi_\al|,
\]
which yields $|y_\al - \xi_\al| \leq C \al\sqrt{\ep}$. Then, using $\Psi(y_\al, \xi_\al,s_\al) \leq \Psi(\hat{x}, \hat x,s_\al)$, \eqref{u-bdd}, and \eqref{v-p-Lip}, we obtain
\begin{align*}
\frac{|\hat x-y_\al|^2}{2\sqrt{\ep}} &\leq u(\hat x,s_\al) - u(y_\al ,s_\al) + \ep\left[ v \left(\frac{\hat x}{\ep}, \frac{\hat z - \hat x}{\sqrt{\ep}} \right) - v \left(\frac{\hat x}{\ep}, \frac{\hat z - \xi_\al}{\sqrt{\ep}} \right)\right]+\gamma\frac{|\hat{x}-\hat{y}|^2}{2}\\
&\leq C|\hat x - y_\al| + C\sqrt{\ep}|\hat{x}-\xi_\al| + C\ep  \\ &\leq C|\hat x - y_\al| + C(1+\alpha)\ep,
\end{align*}
which yields $|\hat x - y_\al|\leq C\sqrt{1+\alpha}\sqrt{\ep}$. Finally, $|\hat{t}-s_{\alpha}|\leq C(1+K)\sqrt{\ep}$ follows from $\Psi(y_{\alpha},\xi_{\alpha},s_{\alpha})\leq \Psi(y_{\alpha},\xi_{\alpha},\hat{t})$, \eqref{u-bdd}, and \eqref{Goal Step1}.

\noindent
\textit{Step 3.2:} We now show that, upon passing to a subsequence, there holds
\begin{align}\label{Step3.2}
(y_\al, \xi_\al,s_\al) \to (\hat{y}, \hat{y}, \hat s)\quad\text{as $\alpha\rightarrow 0^+$}.
\end{align}
In view of \eqref{Step3.1}, there exists $(\tilde{y},\tilde{s})\in \R^n \times [0,\infty)$ such that, upon passing to a subsequence, $(y_\al, \xi_\al,s_\al) \to (\tilde{y}, \tilde{y}, \tilde s)$ as $\al \to 0^+$. As $\Psi(\tilde{y},\tilde{y},\tilde{s}) = \lim_{\alpha\rightarrow 0^+} \Psi(y_{\alpha},\xi_{\alpha},s_{\alpha}) \leq \Psi(\hat y,\hat y,\hat s)$ and using that, by (Step 0), the point $(\hat{y},\hat{s})$ is a strict global minimum of the map $(y,s)\mapsto\Psi(y,y,s)$, we find that $(\tilde{y},\tilde{s})= (\hat{y},\hat{s})$.

\noindent
\textit{Step 3.3:} Introducing $\psi:\R^n\times [0,\infty)\rightarrow \R$ defined by $\psi(x,t):=u(x,t)-\Psi(x,\xi_{\alpha}, t)$, we see that $u-\psi$ has a global minimum at the point $(y_{\alpha},s_{\alpha})$. We compute
\begin{align*}
\psi_t(y_{\alpha},s_{\alpha}) = \frac{\hat t - s_{\alpha}}{\sqrt{\ep}}-K+\gamma(\hat s - s_{\alpha}),\;\;
D\psi(y_{\alpha},s_{\alpha}) = \frac{\hat x - y_{\alpha}}{\sqrt{\ep}}+\frac{\xi_{\alpha}-y_{\alpha}}{\alpha}+\gamma(\hat y - y_{\alpha}).
\end{align*}
By the viscosity supersolution test, local Lipschitz continuity of $\ol H$, and \eqref{Step3.1}, we have for any $\alpha\in (0,1)$ that
\begin{align*}
K&- \frac{\hat t - s_{\alpha}}{\sqrt{\ep}} -\gamma(\hat s - s_{\alpha}) - \ol H\left(\frac{\hat x - y_{\alpha}}{\sqrt{\ep}} + \gamma(\hat y - y_{\alpha}) \right) \\ &= - \psi_t(y_{\alpha},s_{\alpha}) -\ol H\left(D\psi(y_{\alpha},s_{\alpha})-\frac{\xi_{\alpha}-y_{\alpha}}{\alpha} \right)\\ 
&\leq \ol H(D\psi(y_{\alpha},s_{\alpha}))-\ol H\left(D\psi(y_{\alpha},s_{\alpha})-\frac{\xi_{\alpha}-y_{\alpha}}{\alpha} \right) \leq C\sqrt{\ep}.
\end{align*}
In view of \eqref{Step3.2}, passing to the limit $\al \to 0^+$ in the above inequality yields \eqref{Goal Step3}.

\medskip 
\noindent \textbf{Step 4:} For the case $\hat t,\hat s>0$, we combine \eqref{Goal Step2} and \eqref{Goal Step3}, use local Lipschitz continuity of $\ol H$, and \eqref{Goal Step1}, to find that
\begin{align*}
2K\leq \ol H\left(\frac{\hat x - \hat y}{\sqrt{\ep}} \right) -\ol H\left(\frac{\hat z - \hat y}{\sqrt{\ep}} \right) + C\sqrt{\ep}+C\gamma\leq C\sqrt{\ep}+C\gamma,
\end{align*}
which is a contradiction if $\gam \leq \frac{1}{2} \sqrt{\ep}$ and $K=K_1 \sqrt{\ep}$ for $K_1>0$ sufficiently large.
Thus, $\hat t=0$ or $\hat s=0$. In either scenario, using the definition of $\Phi$, the fact that $u^{\ep}(\cdot,0)=u(\cdot,0)$, \eqref{u-bdd} and \eqref{v-p-Lip}, we have for any $(x,t)\in \R^n\times [0,T]$ that
\[
\Phi(x,x,x,t,t)\leq\Phi(\hat a) \leq u^\ep(\hat x,\hat t) - u(\hat y,\hat s) - \ep v\left(\frac{\hat x}{\ep},\frac{\hat z-\hat y}{\sqrt{\ep}}\right) \leq C \sqrt{\ep}.
\]
In view of the definition of $\Phi$, letting $\gam\to 0^+$ in the above inequality yields
\begin{align*}
u^\ep(x,t) - u(x,t) - \ep v\left(\frac{ x}{\ep}, 0\right) - 2K_1\sqrt{\ep}t  \leq C \sqrt{\ep}\qquad \forall (x,t)\in \R^n\times [0,T].
\end{align*}
Finally, by \eqref{v-p-Lip}, we conclude that 
\[
u^\ep(x,t) - u(x,t) \leq C (1+T) \sqrt{\ep}\qquad \forall (x,t) \in \R^n\times [0,T].
\]
\end{proof}

To complete the proof of \eqref{error-bound}, it remains to show that
\begin{align}\label{lowerbd}
u^\ep(x,t) - u(x,t) \geq -C (1+T) \sqrt{\ep}\qquad \forall (x,t) \in \R^n \times [0,T].
\end{align}
\begin{proof}[Proof of \eqref{lowerbd}]
For $K>0$ to be chosen and $\gam >0$ small,
we introduce the auxiliary function $\widetilde\Phi_1:\R^n\times \R^n\times \R^n\times [0,\infty)\times [0,\infty)\rightarrow \R$ given by
\begin{align*}
\widetilde \Phi_1(x,y,z,t,s) := u^\ep(x,t) - u(y,s) - \ep v\left(\frac{x}{\ep},\frac{z-y}{\sqrt{\ep}}\right) + \omega(x,y,z,t,s),
\end{align*}
where $\omega$ is defined as in \eqref{omega def}. By following closely and carefully the proof of \eqref{upperbd}, we obtain the desired result.
\end{proof}

\subsection{Part II: Removal of assumption \eqref{u-ep-bdd}}\label{subsec: Removal}
Let $T>0$ be fixed, $\ep\in (0,1)$, and suppose that we are in the situation (A1)--(A3), i.e., (A1)--(A2) hold and $g\in C^{0,1}(\R^n)$. Let $\rho\in C_c^\infty(\R^n,[0,\infty))$ be a standard mollifier, i.e.,
\[
\int_{\R^n}\rho(x)\,\mathrm{d}x=1, \qquad {\rm supp}(\rho) \subset \{x\in \R^n:|x|\leq 1\}, \qquad \rho(x) = \rho(-x)\text{ for }x\in \R^n,
\]
We set $\rho^\ep:=\frac{1}{\ep^n} \rho(\frac{\cdot}{\ep})$ and $g^{\ep}:=\rho^\ep * g$. Then, $g^\ep\in C^2(\R^n)$ and we have the bounds
\begin{equation}\label{eq:cor-nd-1}
\|g^\ep -g\|_{L^\infty(\R^n)} \leq C\ep,\qquad \|Dg^\ep\|_{L^\infty(\R^n)} + \ep \|D^2g^\ep\|_{L^\infty(\R^n)} \leq C.
\end{equation}
Let $\tilde u^\ep$ denote the viscosity solution to
\begin{equation}\label{1.1with gep}
\begin{cases}
 \tilde u_t^\ep+H\left(\frac{x}{\ep},D\tilde u^\ep\right)=\ep \Delta \tilde u^\ep \qquad &\text{in} \ \R^n \times (0,\infty),\\
\tilde u^\ep(x,0)=g^{\ep}(x) \qquad &\text{on} \ \R^n,
\end{cases} 
\end{equation}
and let $\tilde u$ denote the viscosity solution to
\begin{equation*}
\begin{cases}
 \tilde u_t+\ol H\left(D\tilde u\right)=0 \qquad &\text{in} \ \R^n \times (0,\infty),\\
\tilde u(x,0)=g^{\ep}(x) \qquad &\text{on} \ \R^n.
\end{cases} 
\end{equation*}
By \eqref{eq:cor-nd-1} and the comparison principle, we have that
\begin{equation}\label{eq:cor-nd-2}
\|\tilde u^\ep -u^\ep\|_{L^\infty(\R^n\times[0,\infty))}+\|\tilde u -u\|_{L^\infty(\R^n\times[0,\infty))} \leq C\ep.
\end{equation}
On the other hand, in view of \eqref{eq:cor-nd-1} we have $\left|H(\frac{x}{\ep},Dg^\ep(x)) - \ep \Del g^\ep(x)\right| \leq C$ for $x\in \R^n$, which yields that $(x,t)\mapsto g^\ep(x) + Ct$ is a supersolution to \eqref{1.1with gep}, and $(x,t)\mapsto g^\ep(x) - Ct$ is a subsolution to \eqref{1.1with gep}.
By the comparison principle,
\[
g^\ep(x) - Ct \leq \tilde u^\ep(x,t) \leq g^\ep(x) + Ct \qquad \forall(x,t) \in \R^n\times[0,\infty),
\]
which implies that $\|\tilde u^\ep_t(\cdot,0)\|_{L^\infty(\R^n)} \leq C$.
As $\tilde u^\ep_t$ solves a linear parabolic equation, we find that $\|\tilde u^\ep_t\|_{L^\infty(\R^n\times[0,\infty))} \leq C$ by the maximum principle.
Then, by the classical Bernstein method (see, e.g.,  \cite[Chapter 1]{Tran}),
\begin{equation*}
\|\tilde u^{\ep}_t\|_{L^{\infty}(\R^n \times [0,\infty))}+\|D\tilde u^{\ep}\|_{L^{\infty}(\R^n \times [0,\infty))}\leq C.
\end{equation*}
Thus, we can assume that (A4) holds, and the proof from Section \ref{subsec: Pf based on} yields
\begin{equation}\label{eq:cor-nd-3}
\|\tilde u^\ep  -\tilde u\|_{L^\infty(\R^n\times[0,T))} \leq C(1+T)\sqrt{\ep}.
\end{equation}
Finally, combining \eqref{eq:cor-nd-2} and \eqref{eq:cor-nd-3}, we see that 
\[
\|u^\ep-u\|_{L^\infty(\R^n \times [0,T])} \leq C(1+T) \sqrt{\ep},
\]
which concludes the proof. \qed

\section{Optimality of the bound in Theorem \ref{thm:nd}} \label{sec:optimality}

In this section, we consider the vanishing viscosity problem \eqref{eq:F-ep} with particular choices of $F$ and $g$.
For $\ep>0$, let $u^\ep$ denote the viscosity solution to  \eqref{eq:F-ep}, and let $u$ denote the viscosity solution to \eqref{eq:C} with $\ol H = F$.
If $F \in \Lip_{\rm loc} (\R^n)$ and $g \in \Lip(\R^n)$, then $\|Du^\ep\|_{L^\infty(\R^n \times [0,\infty))} \leq \|Dg\|_{L^\infty(\R^n)}$.
Besides, it was obtained in \cite{Fl, CL, Ev2} that
\begin{equation}\label{eq:sqrt-t-ep}
\|u^\ep-u\|_{L^\infty(\R^n \times [0,T])} \leq C\sqrt{T\ep},
\end{equation}
where the constant $C>0$ depends only on $n$, $R_0:=\|g\|_{\Lip(\R^n)}$, and $\|F\|_{\Lip(B(0,R_0))}$. 

\smallskip

Below, we study both Cauchy problems and Dirichlet problems.

\subsection{A linear Cauchy problem}
We first consider the case that $F$ is linear in one dimension.

\begin{prop}\label{prop:vv-l}
Let $n=1$ and assume that $H(y,p) = F(p)=p$ for $y,p\in \R$.
Let $g\in \Lip(\R)$ and suppose that $g\geq g(-1)$ on $\R$ and $g(x) \geq x+1+g(-1)$ for any $x\in [-1,0]$.
Then, for any $\ep \in (0,\frac{1}{4})$ there holds
\[
|u^\ep(0,1)-u(0,1)| = u^\ep(0,1)-g(-1) \geq \frac{\mathrm{e}-1}{\sqrt{\pi} \mathrm{e}} \sqrt{\ep},
\]
where $u^\ep$ denotes the viscosity solution to \eqref{eq:F-ep} and $u$ denotes the viscosity solution to \eqref{eq:C} with $\ol H = F$.
\end{prop}

\begin{proof}
In this situation, the problem \eqref{eq:F-ep} becomes
\begin{equation*}
\begin{cases}
 u_t^\ep+u^\ep_x=\ep  u^\ep_{xx} \qquad &\text{in} \ \R \times (0,\infty),\\
u^\ep(x,0)=g(x) \qquad &\text{on} \ \R,
\end{cases} 
\end{equation*}
and the problem \eqref{eq:C} becomes
\begin{equation}\label{vv-l-0}
\begin{cases}
 u_t+u_x=0 \qquad &\text{in} \ \R \times (0,\infty),\\
u(x,0)=g(x) \qquad &\text{on} \ \R.
\end{cases} 
\end{equation}
Note that the solution to \eqref{vv-l-0} is given by $u(x,t) = g(x-t)$ for $(x,t) \in  \R \times [0,\infty)$. We introduce $v^{\ep}:\R\times [0,\infty)\rightarrow \R$ given by $v^\ep(x,t):=u^\ep(x+t,t)$. Then, $v^\ep$ solves
\begin{equation*}
\begin{cases}
v_t^\ep=\ep  v^\ep_{xx} \qquad &\text{in} \ \R \times (0,\infty),\\
v^\ep(x,0)=g(x) \qquad &\text{on} \ \R,
\end{cases} 
\end{equation*}
and hence, $v^{\ep}$ is given by
\[
v^\ep(x,t)= \frac{1}{\sqrt{4\pi \ep t}} \int_{-\infty}^\infty \mathrm{e}^{-\frac{|x-y|^2}{4 \ep t}} g(y)\,\mathrm{d}y.
\]
This implies that
\[
u^\ep(x,t) = v^{\ep}(x-t,t) =  \frac{1}{\sqrt{4\pi \ep t}} \int_{-\infty}^\infty \mathrm{e}^{-\frac{|x-y|^2}{4 \ep t}} g(y-t)\,\mathrm{d}y.
\]
{Using $\tilde{g}(x):=g(x)-g(-1)\geq 0$ for all $x\in \R$ and $\tilde{g}(y-1) \geq y$ if $y\in [0,1]$, we have for any $\ep\in(0,\frac{1}{4})$ that
\begin{align*}
u^\ep(0,1) - u(0,1) &= \frac{1}{\sqrt{4\pi \ep}} \int_{-\infty}^\infty \mathrm{e}^{-\frac{y^2}{4 \ep}} g(y-1)\,\mathrm{d}y - g(-1) \\
&=\frac{1}{\sqrt{4\pi \ep}} \int_{-\infty}^\infty \mathrm{e}^{-\frac{y^2}{4 \ep}} \tilde{g}(y-1)\,\mathrm{d}y \\ &\geq \frac{1}{\sqrt{4\pi \ep}} \int_{0}^{2\sqrt{\ep}} \mathrm{e}^{-\frac{y^2}{4 \ep}} y\,\mathrm{d}y = \frac{\mathrm{e}-1}{\sqrt{\pi} \mathrm{e}}  \sqrt{\ep}. 
\end{align*}
Noting that $u(0,1) = g(-1)$ yields the desired result.}
\end{proof}

\subsection{Proof of Theorem \ref{thm:new}}
Recall that we consider \eqref{eq:F-ep} in one dimension, i.e.,
\begin{equation}\label{eq:new-ep}
\begin{cases}
 u_t^\ep+F\left(u^\ep_x\right)=\ep  u^\ep_{xx} \qquad &\text{in} \ \R \times (0,\infty),\\
u^\ep(x,0)=g(x) \qquad &\text{on} \ \R.
\end{cases} 
\end{equation}
As $\ep \to 0^+$, we have that $u^\ep \to u$ locally uniformly on $\R \times [0,\infty)$, where $u$ solves
\begin{equation}\label{eq:new}
\begin{cases}
 u_t+F\left(u_x\right)=0 \qquad &\text{in} \ \R \times (0,\infty),\\
u(x,0)=g(x) \qquad &\text{on} \ \R.
\end{cases} 
\end{equation}

\begin{proof}[Proof of Theorem \ref{thm:new}]
We first show that $u(0,1)=0$.
As $F(0) = 0$ and $g\geq 0$ on $\R$, the function $\varphi  \equiv 0$ is a subsolution to \eqref{eq:new}, which yields $u(0,1) \geq 0$. In order to show that also $u(0,1)\leq 0$, let us introduce $h,\tilde{h}:\R\rightarrow \R$ given by $\tilde{h}(x):=x+1$ and $h(x) := \max\{\tilde{h}(x), 0\}$ for $x \in \R$. Let $\rho \in C_c^\infty(\R, [0,\infty))$ be a standard mollifier, i.e.,
\[
\int_{-\infty}^\infty \rho(x)\,\mathrm{d}x=1,  \quad {\rm supp}(\rho) \subset [-1,1], \quad \text{ and } \quad \rho(x)=\rho(-x) \text{ for } x\in \R.
\]
For $\delta \in (0,1)$, let $\rho^\delta := \frac{1}{\delta} \rho(\frac{\cdot}{\delta})$ and set $h^\delta := \rho^\delta * h$. Note that $h^\delta \geq 0$ as $h\geq 0$, and $h^\delta \geq \rho^{\delta}*\tilde{h} = \tilde{h}$ as $h\geq \tilde{h}$. Hence, $h^\delta \geq h \geq g$ on $\R$. Besides, $0 \leq (h^\delta)' \leq 1$ on $\R$ which follows from $(h^\delta)' = \rho^{\delta} * h'$ and $0\leq h'\leq 1$ a.e. on $\R$. Introducing $\psi(x,t) := h(x-t)$ and $\psi^\delta(x,t) := h^\delta(x-t)$ for $(x,t) \in \R \times [0,\infty)$, we have that $\psi^\delta(\cdot,0) = h^{\delta}\geq g$ on $\R$ and, using that $F(p)=p$ for $p\in [0,1]$,
\begin{align*}
\psi^\delta_t(x,t) + F(\psi^\delta_x(x,t))&= -  (h^\delta)'(x-t) + F((h^\delta)'(x-t))=0,
\end{align*}
i.e., $\psi^\delta$ is a supersolution to \eqref{eq:new}. Since $\psi^\delta  \to \psi$ locally uniformly on $\R \times [0,\infty)$ as $\delta \to 0^+$, we deduce that $\psi$ is also a supersolution to \eqref{eq:new}.
In particular, 
\[
u(0,1) \leq \psi(0,1) = h(-1)=0.
\]
Thus, $u(0,1)=0$.

Next, we construct a subsolution to \eqref{eq:new-ep}. We set
\begin{align}\label{phieps def}
\phi^\ep(x,t) :=  \frac{1}{\sqrt{4\pi \ep t}} \int_{-\infty}^\infty \mathrm{e}^{-\frac{|x-y|^2}{4 \ep t}} g(y-t)\,\mathrm{d}y
\end{align}
and recall from the proof of Proposition \ref{prop:vv-l} that $\phi^\ep_t + \phi^\ep_x = \ep \phi^\ep_{xx}$ in $\R\times (0,\infty)$ and $\phi^{\ep}(\cdot,0) = g$ on $\R$. Since $|g'|\leq 1$ a.e. on $\R$, we note that $|\phi^\ep_x|\leq 1$ in $\R\times (0,\infty)$. Using that by assumption $F(p) \leq p$ if $|p|\leq 1$, we find that
\[
\phi^\ep_t + F(\phi^\ep_x) - \ep \phi^\ep_{xx} \leq \phi^\ep_t + \phi^\ep_x - \ep \phi^\ep_{xx} = 0\quad\text{in }\R\times (0,\infty).
\]
Therefore, $\phi^\ep$ is a subsolution to \eqref{eq:new-ep} and by the comparison principle we have that $u^\ep \geq \phi^\ep$. Hence, for $\ep\in (0,\frac{1}{4})$ there holds
\[
u^\ep(0,1) \geq  \phi^\ep(0,1) \geq \frac{\mathrm{e}-1}{\sqrt{\pi} \mathrm{e}}  \sqrt{\ep},
\]
where the second inequality follows from Proposition \ref{prop:vv-l} applied to $\phi^\ep$.
\end{proof}

\begin{rem}\label{rem:C2}
Fix {$\ep\in (0,\frac{1}{5})$} and $\alpha \in (0,\frac{1}{10}\sqrt{\ep})$. In the situation of Theorem \ref{thm:new}, if we replace the initial condition $g$ by
\[
g^\alpha := \rho^\alpha*g \in C^\infty(\R),
\]
where $\rho^\alpha:=\frac{1}{\alpha}\rho(\frac{\cdot}{\alpha})$ with $\rho$ as in the proof of Theorem \ref{thm:new}, then we still have that 
\[
u^\ep(0,1) - u(0,1) \geq \frac{\mathrm{e}-1}{2\sqrt{\pi} \mathrm{e}}  \sqrt{\ep}. 
\]
Indeed, as $g^\alpha(-1) \in (0, \alpha)$ and $g^\alpha(x) \geq 1+x$ for all $x\in [-1,-1+2\sqrt{\ep} ]\subset [-1,-\al)$, we have in view of the proof of Proposition \ref{prop:vv-l} that $u^\ep(0,1) \geq \frac{\mathrm{e}-1}{\sqrt{\pi} \mathrm{e}}  \sqrt{\ep}$ and that $u(0,1) = g^\alpha(-1) <  \frac{1}{10}\sqrt{\ep}$. Note that we still have $g^\al\in \Lip(\R)$ with $\|g^\alpha\|_{L^\infty(\R)} \leq 1$ and $\|(g^\alpha)'\|_{L^\infty(\R)} \leq 1$.
\end{rem}

We now provide a generalization of Theorem \ref{thm:new}.

\begin{prop}\label{prop:generalization}
Let $n=1$ and $\ep \in (0,\frac{1}{4})$.
Let $m>1$ be a constant such that
\[
1-\frac{1}{m} \leq  \frac{e-1}{2\sqrt{\pi} e} \sqrt{\ep}.
\]
Let $F \in \Lip_{\rm loc}(\R)$ be such that 
\[
\begin{cases}
F(p)=\frac{1}{m} p^m \qquad &\text{ for $p\in [0,1]$}, \\
F(p) \leq p \qquad &\text{ for $p\in [-1,0]$}.
\end{cases}
\]
Assume that $g(x)=\max\{1-|x|,0\}$ for $x\in \R$.
Let $u^\ep$ denote the viscosity solution to  \eqref{eq:F-ep} and let $u$ denote the viscosity solution to \eqref{eq:C} with $\ol H = F$.
Then,
\[
|u^\ep(0,1)-u(0,1)|  \geq \frac{e-1}{2\sqrt{\pi} e} \sqrt{\ep}.
\]

\end{prop}

\begin{proof}
First, we note that $F(p) \leq p$ if $|p|\leq 1$.
Thus, by the last part of the proof of Theorem \ref{thm:new}, we still have $u^\ep \geq \phi^\ep$, and hence,
\begin{align}\label{uepslowerbd}
u^\ep(0,1) \geq  \phi^\ep(0,1) \geq \frac{e-1}{\sqrt{\pi} e}  \sqrt{\ep},
\end{align}
where $\phi^\ep$ is defined in \eqref{phieps def}. Now, let $h^\delta = \rho^\delta*h$ for $\delta \in (0,1)$ be defined as in the proof of Theorem \ref{thm:new}.
We recall that $h^\delta \geq h \geq g$ and $0 \leq (h^\delta)' \leq 1$ on $\R$. Introducing $\zeta(x,t) := h(x-t) + (1-\frac{1}{m})t$ and
\[
\zeta^\delta(x,t) := h^\delta(x-t) +\left(1-\frac{1}{m}\right)t \quad \text{ for } (x,t) \in \R \times [0,\infty),
\]
we claim that $\zeta^\delta$ is a supersolution to \eqref{eq:new}.
Indeed, we have that $\zeta^\delta(\cdot,0) = h^\delta \geq g$ on $\R$, and using $F(p)=\frac{1}{m} p^m$ for $p\in [0,1]$ and Bernoulli's inequality, there holds
\begin{align*}
\zeta^\delta_t(x,t) + F(\zeta^\delta_x(x,t))&= 1-\frac{1}{m}-  (h^\delta)'(x-t) + F((h^\delta)'(x-t))\\
&=\frac{((h^\delta)'(x-t))^m-[1+m( (h^\delta)'(x-t)-1)]}{m}\geq 0.
\end{align*}
Since $\zeta^\delta  \to \zeta$ locally uniformly on $\R \times [0,\infty)$ as $\delta \to 0^+$, we deduce that $\zeta$ is also a supersolution to \eqref{eq:new}.
Hence, we have that 
\[
u(0,1) \leq \zeta(0,1) = h(-1)+1-\frac{1}{m} =1-\frac{1}{m}  \leq  \frac{\mathrm{e}-1}{2\sqrt{\pi} \mathrm{e}} \sqrt{\ep},
\]
which completes the proof in view of \eqref{uepslowerbd}.
\end{proof}
It is important to note that in the above proposition, although $F$ depends on $m$ and hence $\ep$, the value $\|F\|_{\Lip_{\rm loc}(\R)}$ does not depend on $m$ and $\ep$.

\subsection{Quadratic Hamiltonian} \label{subsec:quadHam}

Next, we consider the case where $F$ is quadratic in one dimension.  First, we construct an example that complements \eqref{eq:sqrt-t-ep} when $T$ is very small. 

\begin{prop}\label{prop:quad-optimal}
Let $n=1$ and assume that $H(y,p)=F(p)=\frac{1}{2} |p|^2$ for $y,p\in \R$, and  $g(x)=-|x|$ for $x\in \R$.
For $\ep>0$, let $u^\ep$ be the viscosity solution to  \eqref{eq:F-ep} and let $u$ be the viscosity solution to \eqref{eq:C}.
Then, for $\ep \in (0,1)$ and $T>0$, there holds
\[
\|u^\ep - u\|_{L^\infty (\R \times [0,T])} \leq C \sqrt{T\ep},
\]
and this upper bound  $O(\sqrt{T\ep})$ is sharp in the sense that 
\[
\lim_{t\to 0^+}\frac{u^{\ep}(0,t)-u(0,t)}{\sqrt{t\ep}}=-\frac{2}{\sqrt{\pi}}.
\]
\end{prop}
It is important to note that we also obtain a rigorous asymptotic expansion of $u^\ep(0,t)$ for $0<t \ll \ep$ in the proof of Proposition \ref{prop:quad-optimal}. We then give a finer bound of $u^\ep-u$ in Proposition \ref{prop:eplog} under some appropriate conditions on $g$.

In this subsection, we assume the setting of Proposition \ref{prop:quad-optimal}.
Then, the problem \eqref{eq:F-ep} reads 
\begin{equation*}
\begin{cases}
u^\ep_t + \frac{1}{2}|u^\ep_x|^2 = \ep u^\ep_{xx} \quad &\text{ in } \R \times (0,\infty),\\
u^\ep(x,0) = g(x) \quad &\text{ on } \R.
\end{cases}
\end{equation*}
We have that $u^\ep \to u$ locally uniformly on $\R \times [0,\infty)$ as $\ep \to 0^+$, and $u$ solves
\begin{equation}\label{eq:example-0}
\begin{cases}
u_t + \frac{1}{2}|u_x|^2 = 0 \quad &\text{ in } \R \times (0,\infty),\\
u(x,0) = g(x) \quad &\text{ on } \R.
\end{cases}
\end{equation}

\begin{proof}[Proof of Proposition \ref{prop:quad-optimal}]
The bound $O(\sqrt{T\ep})$ was obtained in \cite{Fl, Ev2}.
We only need to show that this bound is optimal here.
As $g(x)=-|x|$ for $x\in \R$, we see that the solution to \eqref{eq:example-0} is given by
\[
u(x,t) = -|x| -\frac{t}{2} \quad \text{ for } (x,t) \in \R \times [0,\infty).
\]
In particular, $u(0,t)=-\frac{t}{2}$ for all $t\geq 0$. For $\ep\in (0,1)$, we have the following representation formula for $u^\ep$ (see, e.g.,  \cite[Chapter 4]{Ev3})
\begin{align}\label{repform}
u^\ep(x,t) = -2\ep \log \left[ \frac{1}{\sqrt{4\pi \ep t}} \int_{-\infty}^\infty \mathrm{e}^{-\frac{|x-y|^2}{4\ep t}-\frac{g(y)}{2\ep}}\,\mathrm{d}y \right].
\end{align}
In particular, for any $t>0$ we have that
\begin{align*}
u^\ep(0,t)&=-2\ep \log \left[ \frac{1}{\sqrt{\pi }} \int_{-\infty}^\infty \mathrm{e}^{-|z|^2+\frac{|z| \sqrt{t}}{\sqrt{\ep}}}\,\mathrm{d}z \right]\\
&=-2\ep \log \left[ \frac{2}{\sqrt{\pi }}\mathrm{e}^{\frac{t}{4\ep}} \int_{0}^\infty \mathrm{e}^{-(z-\frac{\sqrt{t}}{2\sqrt{\ep}})^2}\,\mathrm{d}z\right]\\
&=-\frac{t}{2} -2\ep \log \left[1+ \frac{2}{\sqrt{\pi }} \int_{0}^{\frac{\sqrt{t}}{2\sqrt{\ep}}} \mathrm{e}^{-s^2}\,\mathrm{d}s\right]
=-\frac{t}{2} -2\ep \log  \left[1 + \erf \left(\frac{\sqrt{t}}{2\sqrt{\ep}}\right) \right].
\end{align*}
Note that we have
\[
\erf(z) = \frac{2}{\sqrt{\pi}} \sum_{k=0}^\infty \frac{(-1)^k z^{2k+1}}{k! (2k+1)}.
\]
For $z=\frac{\sqrt{t}}{2\sqrt{\ep}} \ll 1$, we have that $ 0< \erf(z) \ll 1$ and
\[
\log ( 1 + \erf(z)) = \log \left ( 1 + \frac{2}{\sqrt{\pi}} z + \cdots \right) = \frac{2}{\sqrt{\pi}} z +\cdots,
\]
and thus,
\[
u^\ep(0,t) = -\frac{t}{2} -2\ep \frac{2}{\sqrt{\pi}}\frac{\sqrt{t}}{2\sqrt{\ep}} + \cdots
= -\frac{t}{2} - \frac{2\sqrt{t\ep}}{\sqrt{\pi}} + \cdots,
\]
which gives us the desired result.
\end{proof}

\subsubsection{Improvement of convergence rates.} It was shown in \cite{TT} that if the  Hamiltonian $F$ is uniformly convex, i.e., $F''\geq \alpha$ on $\R$ for some constant $\alpha > 0$, then the convergence rate of the vanishing viscosity limit can be improved to $O(\ep |\log \ep|)$ when the initial datum satisfies certain technical assumptions.  
Below  we show that for $F(p) = \frac{1}{2}|p|^2$,  the rate is almost everywhere $O(\ep)$ when the initial datum is $C^2$, although $O(\ep|\log \ep|)$ could happen at some points. For $F(p) = \frac{1}{2}|p|^2$, we have by the Hopf-Lax formula (see e.g., \cite{Ev3, Tran}) that
\begin{align*}
u(x,t)= \inf_{y\in \R} \left\{g(y)+t\, L\left(\frac{x-y}{t}\right)\right\},\quad\text{where}\quad L(v):=\sup_{p\in \R} \left\{pv-F(p)\right\} = \frac{1}{2}|v|^{2}.
\end{align*}
In particular, for any $(x,t)\in \R \times (0,\infty)$ there holds
\begin{align}\label{formulauxt}
u(x,t) - \frac{|x|^2}{2t} =  \inf_{y\in \R} \left\{g(y)+\frac{|y|^2}{2t}-\frac{xy}{t}\right\},
\end{align}
which yields that $x \mapsto u(x,t) - \frac{|x|^2}{2t} $ is concave for any fixed $t>0$. Hence, for any fixed $t>0$, we have that $u(\cdot,t)$ is twice differentiable a.e. on $\R$. For $t>0$, we set
\begin{align}\label{St defn}
S_t:=\left\{x\in \R\,:\, \text{$u(\cdot,t)$ is twice differentiable at $x$}\right\},
\end{align}
and note that $\R\backslash S_t$ has Lebesgue measure zero. Assume now that $g\in C^2(\R)$.
We know that for each $x\in S_t$,   there exists a unique $y_{x,t}\in \R$ such that
\begin{align}\label{yxt}
u(x,t)=\inf_{y\in \R}\left\{g(y)+\frac{1}{2t}|x-y|^2\right\} = g(y_{x,t}) + \frac{1}{2t} |x-y_{x,t}|^2,
\end{align}
and we have that
\begin{align}\label{yxt prop}
u_x(x, t)=g'(y_{x,t})=\frac{x-y_{x,t}}{t},\qquad g''(y_{x,t})\geq -\frac{1}{t},
\end{align}
where the first equality in \eqref{yxt prop} follows from the method of characteristics. We obtain that $u_x(x, t) = g'( x-tu_x(x, t))$ and hence,
\[
u_{xx}(x,t) = \left( 1 - t u_{xx}(x,t) \right)g''(x-t u_x(x,t))  = \left( 1 - t u_{xx}(x,t) \right)g''(y_{x,t}).
\]
In view of \eqref{yxt prop}, we deduce that 
\begin{align}\label{yxt prop2}
g''(y_{x,t}) > -\frac{1}{t}.
\end{align}
We are now in a position to prove the following result:

\begin{prop}\label{prop:eplog}
Let $n=1$ and assume that $H(y,p)=F(p)=\frac{1}{2} |p|^2$ for $y,p\in \R$.  Assume that $g \in \mathrm{Lip}(\R)$. For $\ep\in (0,1)$, let $u^\ep$ denote the viscosity solution to  \eqref{eq:F-ep} and let $u$ denote the viscosity solution to \eqref{eq:C}  with $\ol H = F$.  Then, the following assertions hold true. 
\begin{itemize}
\item[(i)] For fixed $(x,t)\in \R\times [0,\infty)$,  there holds
\begin{equation*}
|u^\ep(x,t) - u(x,t)|\leq 2\ep |\log \ep|
\end{equation*}
for $\ep>0$ sufficiently small.
\item[(ii)] If we further assume that $g\in C^2(\R)$, then for each fixed $t>0$ we have that 
\begin{equation*}
|u^\ep(x,t) - u(x,t)|\leq C \ep \quad \text{for a.e. $x\in \R$},
\end{equation*}
where $C=C(x,t)>0$ is independent of $\ep \in (0,1)$.

\item[(iii)] If $g(x)=-\frac{1}{2}x^2$ for all $x\in [-1,1]$ and $g(x) \geq -\frac{1}{2}x^2$ for all $x\in \R$, then 
\begin{equation*}
|u^\ep(0,1) - u(0,1)|\geq \frac{1}{2}\ep |\log \ep| 
\end{equation*}
for $\ep>0$ sufficiently small.
\end{itemize}
\end{prop}

\begin{proof} 
Without loss of generality, let $(x,t)=(0,1)$. Introducing $h(y):=g(y)+\frac{1}{2}|y|^2$ for $y\in \R$, we have by \eqref{repform} and \eqref{formulauxt} that
\begin{align*}
u^\ep(0,1) = -2\ep \log \left[ \frac{1}{\sqrt{4\pi \ep}} \int_{-\infty}^\infty \mathrm{e}^{-\frac{h(y)}{2\ep}}\,\mathrm{d}y \right],\qquad u(0,1)=\min_{\R}h = h(\bar y),
\end{align*}
where $\bar y\in \R$ is a fixed point for which there holds $h(\bar y)=\min_{\R} h$. Note that
\begin{align}\label{diffformula}
u^\ep(0,1) - u(0,1) =-2\ep \log \left[ \frac{1}{\sqrt{4\pi \ep }} \int_{-\infty}^\infty \mathrm{e}^{-\frac{h(y)-h(\bar y)}{2\ep}}\,\mathrm{d}y \right] .
\end{align}
We first prove (i).  
Since $g\in \mathrm{Lip}(\R)$, there exists $M>0$ such that for any $y\in \R$ with $|y-\bar y|\geq M$ there holds $h(y)-h(\bar y)\geq \frac{1}{4}|y-\bar y|^2$. For $\ep \in (0,1)$, we have that
\[
2M \geq \int_{\bar y -M}^{\bar y+M} \mathrm{e}^{-\frac{h(y)-h(\bar y)}{2\ep}}\,\mathrm{d}y \geq \int_{\bar y-M}^{\bar y+M} \mathrm{e}^{-\frac{A|y-\bar y|}{2\ep}}\,\mathrm{d}y = \frac{4}{A}\left(1-\mathrm{e}^{-\frac{AM}{2\ep}}\right)\ep \geq B\ep
\]
for some $A=A(L, |\bar y|, M)>0$   and $B=B(A,M)>0$.  Moreover, 
\begin{align*}
 \int_{\R\backslash (\bar{y}-M,\bar{y}+M)} \mathrm{e}^{-\frac{h(y)-h(\bar y)}{2\ep}}\,\mathrm{d}y  \leq  \int_{\R\backslash (\bar{y}-M,\bar{y}+M)} \mathrm{e}^{-\frac{|y-\bar y|^2}{8\ep}}\,\mathrm{d}y \leq  2\int_{M}^\infty \mathrm{e}^{-\frac{M y}{8\ep}}\,\mathrm{d}y   \leq  \frac{ 16}{M}\ep.
\end{align*}
Combining the two inequalities stated above, we find that
\[
B\ep\leq \int_{-\infty}^\infty \mathrm{e}^{-\frac{h(y)-h(\bar y)}{2\ep}}\,\mathrm{d}y \leq 2M+ \frac{ 16}{M}\ep.
\]
Thus, in view of \eqref{diffformula}, we obtain that
\[
 |u^\ep(0,1) - u(0,1)| \leq 2 \ep |\log \ep|
\]
for $\ep > 0$ sufficiently small.

Next we prove (ii).   Assume that  $0\in S_1$, where $S_1\subset \R$ is defined in \eqref{St defn}. Then, $\bar y=y_{0,1}$ (recall \eqref{yxt}--\eqref{yxt prop}) is the unique minimum point of $h$ and, using \eqref{yxt prop2}, we have that $h''(\bar y)=g''(\bar y)+1>0$. Combining with the fact that $g\in C^2(\R) \cap \Lip(\R)$,  there exists $\alpha>0$ such that $\alpha |y-\bar y|^2 \leq h(y)-h(\bar y)\leq \frac{1}{\alpha} |y-\bar y|^2$ for any $y\in \R$. Thus, there exists $C=C(\alpha)>0$ such that
\[
 {C}\sqrt{\ep}\leq  \int_{-\infty}^\infty \mathrm{e}^{-\frac{h(y)-h(\bar y)}{2\ep}}\,\mathrm{d}y \leq  \frac{1}{C}\sqrt{\ep},
\]
and hence, in view of \eqref{diffformula}, 
\[
|u^\ep(0,1) - u(0,1)|\leq C \ep.
\]

Finally, (iii) follows immediately from \eqref{diffformula} in combination with the fact that due to $h \geq 0$ on $\R$, $\left. h\right\rvert_{[-1,1]}\equiv 0$, and $h(\bar y)=0$, there holds
\[
\int_{-\infty}^\infty \mathrm{e}^{-\frac{h(y)-h(\bar y)}{2\ep}}\,\mathrm{d}y \geq \int_{-1}^1 \mathrm{e}^{-\frac{h(y)-h(\bar y)}{2\ep}}\,\mathrm{d}y=2. 
\]
\end{proof}

It is not clear to us whether Proposition  \ref{prop:eplog} holds for other uniformly convex $F$.  For strictly but not uniformly convex $F$ (e.g.,  $F(p)=\frac{1}{4}|p|^4$), the convergence rate for the vanishing viscosity process might be $O(\ep^{\alpha})$ for some exponent $\alpha\in (\frac{1}{2}, 1)$; see the numerical Example \ref{ex:5}.  A natural question is whether we will see a similar convergence rate when the homogenization process is involved for the quadratic case as numerical Example \ref{ex:10} suggests.  Let us briefly demonstrate the technical difficulty in extending the proof of Proposition  \ref{prop:eplog} to the homogenization problem. Consider  $H(y,p)=\frac{1}{2}|p|^2+V(y)$ for a smooth  $\Z^n$-periodic potential function $V$.  Then, by the Hopf-Cole transformation, we have that
\[
u^{\ep}(x,t)=-2\ep \log \left[ h\left(\frac{x}{\ep}, \frac{t}{\ep}\right)\right],
\]
where $h = h(x,t)$ is the solution to the problem
\begin{equation*}
\begin{cases}
h_t-\Delta h+\frac{1}{2}Vh=0 \qquad &\text{in} \ \R^n \times (0,\infty),\\
h(x,0)=\mathrm{e}^{-\frac{g(\ep x)}{2\ep}} \qquad &\text{on} \ \R^n.
\end{cases} 
\end{equation*} 
Therefore, we have that
\[
u^{\ep}(x,t)=-2\ep \log\left[\int_{\R^n}K\left(\frac{x}{\ep}, \frac{y}{\ep}, \frac{t}{\ep}\right)\mathrm{e}^{-\frac{g(y)}{2\ep}}\,\mathrm{d}y\right],
\]
where $K = K(x,y,t)$ denotes the fundamental solution corresponding to the operator $\partial_t-\Delta +\frac{1}{2}V$. Obtaining the convergence rate requires a sharp estimate of the homogenization of $K$, which is a highly nontrivial subject.  
%When $V\in L^r(\R^n)$ for some $r>0$ or when $V$ satisfies certain decay conditions,  there are fruitful results in comparing $K(x,y,t)$ with the regular Gaussian kernel $(4\pi t)^{-\frac{n}{2}}\mathrm{e}^{-\frac{|y-x|^2}{4t}}$; see \cite{BDS}  and the references therein. For periodic $V$, the only estimate we are aware of is the Li-Yau inequality \cite{LY} for the  operator $\partial_t-\Delta +\frac{1}{2}V$ defined on a manifold, which is however not delicate enough for our purpose. 
 Let us further point out that the convergence rate might also depend on the regularity of $V$ as the numerical Example \ref{ex:6} suggests.

\subsection{A Dirichlet problem}
We are again in one dimension.
For $\ep > 0$, we consider the Dirichlet problem
\begin{equation}\label{D-ep}
\begin{cases}
2 (u^\ep)^3 = \ep (u^\ep)'' \quad \text{ in } (0,\infty),\\
u^\ep(0)=1.
\end{cases}
\end{equation}
It is quickly seen that the solution is given by
\[
u^\ep(x) = \frac{\sqrt{\ep}}{x+\sqrt{\ep}} \quad \text{ for } x \geq 0.
\]
In particular, we have $u^\ep \to u \equiv 0$ locally uniformly in $(0,\infty)$.
Of course, there is a boundary layer of size $O(\sqrt{\ep})$ at $x=0$, but let us ignore this boundary layer in our discussion here.
We observe that for any $\ep \in (0,1)$ there holds
\[
|u^\ep(1)-u(1)| =  \frac{\sqrt{\ep}}{1+\sqrt{\ep}} \geq \frac{1}{2}\sqrt{\ep}.
\]
Thus, once again, we see that the $O(\sqrt{\ep})$ rate occurs naturally here.
We record this in the following lemma.
\begin{lem}\label{lem:D}
For $\ep \in (0,1)$, let $u^\ep$ denote the solution to \eqref{D-ep}.
Then, $u^\ep \to u \equiv 0$ locally uniformly in $(0,\infty)$, and there holds
\[
|u^\ep(x) - u(x)| =  \frac{\sqrt{\ep}}{x+\sqrt{\ep}} \geq \frac{1}{2x}\sqrt{\ep}\qquad \forall x\geq 1.
\]
In particular, for any $d>1$, the optimal rate for the convergence of $u^\ep$ to $u$ in the $L^\infty((1,d))$-norm is $O(\sqrt{\ep})$.
\end{lem}

\section{Numerical results for the vanishing viscosity process and the homogenization problem} \label{sec:num}

\subsection{Vanishing viscosity process}
We consider \eqref{eq:F-ep} in one dimension, that is,
\begin{equation}\label{vv}
\begin{cases}
 u_t^\ep+F\left(u^\ep_x\right)=\ep  u^\ep_{xx} \qquad &\text{in} \ \R \times (0,\infty),\\
u^\ep(x,0)=g(x) \qquad &\text{on} \ \R.
\end{cases} 
\end{equation}
Recall that, as $\ep \to 0^+$, $u^\ep \to u$ locally uniformly on $\R \times [0,\infty)$, where $u$ solves
\begin{equation*}
\begin{cases}
 u_t+F\left(u_x\right)=0 \qquad &\text{in} \ \R \times (0,\infty),\\
u(x,0)=g(x) \qquad &\text{on} \ \R.
\end{cases} 
\end{equation*}
We now verify numerically that, in some particular examples,
\[
\|u^\ep(\cdot,1)- u(\cdot,1)\|_{L^\infty} \geq C \sqrt{\ep},
\]
for some $C>0$ independent of $\ep \in (0,1)$, which  confirms again that the bound $O(\sqrt{\ep})$ is optimal in general.
To do so, we consider various choices of $F$ and $g$ and compute $\|u^\ep(\cdot,1)- u(\cdot,1)\|_{L^\infty}$ for different values of $\ep>0$.

Let us describe our methodology.
We partition a spatial interval $[a, b]$ by a uniform mesh with mesh size $\Delta{x}$ and choose adaptive time steps $\Delta{t}$ to march through a given time interval $[0,T]$. 
Accordingly, we discretize equation \eqref{vv} as follows:
\begin{align*}
u_i^{n+1} = u_i^n -\Delta{t}\left[F\left(\frac{u_{i+1}^n-u_{i-1}^n}{2\Delta{x}}\right)-\ep\frac{u_{i+1}^n-2u_i^n+u_{i-1}^n}{\Delta{x}^2} \right] =:G(u_{i-1}^n,u_i^n,u_{i+1}^n).
\end{align*}
Monotonicity of the scheme requires that $G$ is nondecreasing in each of its arguments; consequently, we have 
\begin{eqnarray}
\ep &\geq& \frac{1}{2}\Delta{x}\max_{p} |F'(p)|,  \label{minvis}\\
\Delta{t} &\leq& \frac{\Delta{x}^2}{2\ep}. \label{timestep} 
\end{eqnarray}
The condition \eqref{minvis} requires a minimum viscosity to be imposed on the numerical scheme, and the time step has to be chosen according to \eqref{timestep}. To check the effect of vanishing viscosity, we will set 
\begin{eqnarray*}
\ep_{\min} &=& \frac{1}{2}\Delta{x}\max_{p}|F'(p)|,  \\
\ep &=&2^{k}\ep_{\min} \quad \mbox{ for } k= 9, \cdots, 1, 0,  \\
\Delta{t} &=& c_{\rm cfl}\frac{\Delta{x}^2}{2\ep},  
\end{eqnarray*}
where $c_{\rm cfl}\leq 1$ is the CFL number. 
We note that it is extremely hard to verify rigorously the examples considered below.

\begin{figure}
	\begin{subfigure}{0.49\textwidth}
		\includegraphics[width=\textwidth]{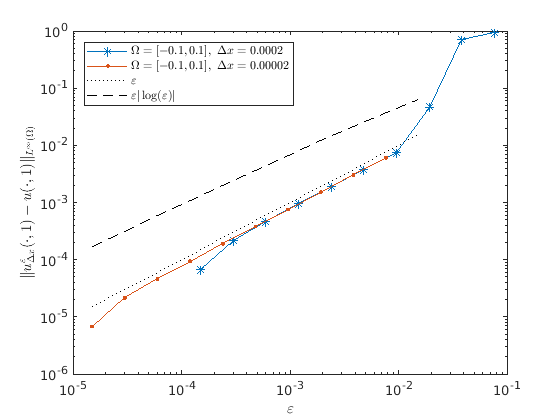}
		\subcaption{$F(p) = |p|^\frac{3}{2}$ for $p\in \R$,\\ $g(x) = -|x|$ for $x\in \R$.}
	\end{subfigure}
	\begin{subfigure}{0.49\textwidth}
		\includegraphics[width=\textwidth]{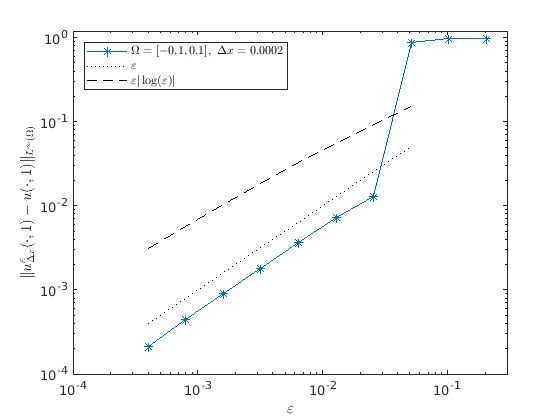}
		\subcaption{$F(p) = |p|^4$ for $p\in \R$,\\ $g(x) = -|x|$ for $x\in \R$.}
	\end{subfigure}
	
\begin{subfigure}{0.49\textwidth}
		\includegraphics[width=\textwidth]{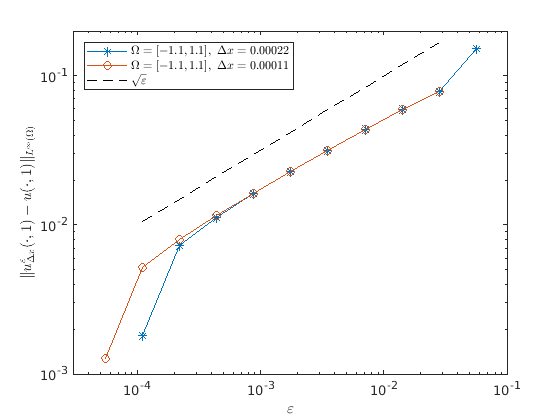}
		\subcaption{$F(p) = |p|$ for $p\in \R$, \\ $g(x) = \max\{1-|x|,0\}$ for $x\in \R$.}
	\end{subfigure}
	\begin{subfigure}{0.49\textwidth}
		\includegraphics[width=\textwidth]{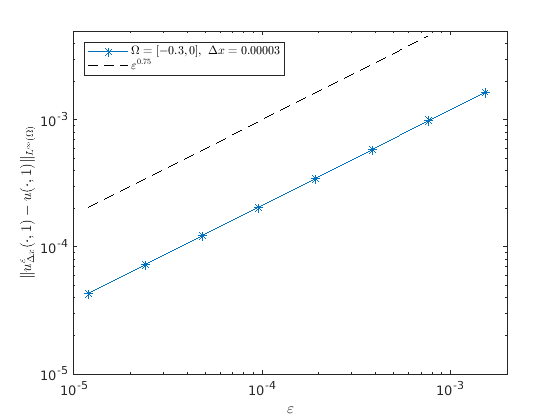}
		\subcaption{$F(p) = p^3$ for $p\in \R$, \\ $g(x) = \frac{2\sqrt{2}}{9} (-x)^{3/2}$ for $x\in (-\infty,0]$.}
	\end{subfigure}	
	
\begin{subfigure}{0.49\textwidth}
		\includegraphics[width=\textwidth]{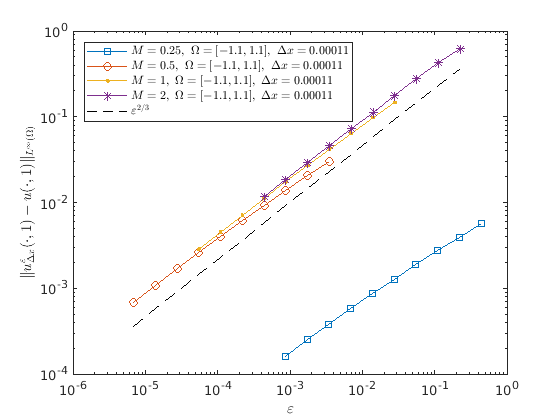}
		\subcaption{$F(p) = \frac{1}{4} |p|^4$ for $p\in \R$, \\ $g(x) = M \min\{|x|,|x-\frac{1}{2}|-\frac{1}{4}\}$ for $x\in \R$.}
	\end{subfigure}
	
	\caption{Illustration of the error $\|u^{\ep}_{\Delta x}(\cdot,1)-u(\cdot,1)\|_{L^\infty(\Omega)}$ for Examples \ref{ex:1}--\ref{ex:5}.}
	\label{Fig: plots}
\end{figure}

\begin{ex} \label{ex:1}
Assume $F(p)=|p|^{\nicefrac{3}{2}}$ for $p\in \R$, and $g(x)=-|x|$ for $x\in \R$. 
Then,
\[
u(x,t)=-|x|-t \quad \text{ for all } (x,t)\in \R\times [0,\infty).
\]
Numerical results are shown in Figure \ref{Fig: plots} (A). We observe that the convergence rate is $O(\ep)$ in this example.
\end{ex}

\begin{ex} \label{ex:2}
Assume $F(p)=|p|^{4}$ for $p\in \R$, and $g(x)=-|x|$ for $x \in \R$. 
Then,
\[
u(x,t)=-|x|-t \quad \text{ for all } (x,t)\in \R\times [0,\infty).
\]
Numerical results are shown in Figure \ref{Fig: plots} (B). We observe that the convergence rate is $O(\ep)$ in this example.
\end{ex}

\begin{ex}\label{ex:3}
Assume $F(p)=|p|$ for $p\in \R$, and $g(x)=\max\{1-|x|,0\}$ for $x\in \R$. 
Then,
\[
u(x,t)=\max\{1-|x|-t, 0\} \quad \text{ for all } (x,t)\in \R\times [0,\infty).
\]
Numerical results are shown in Figure \ref{Fig: plots} (C). We observe that the convergence rate is $O(\sqrt{\ep})$ in this example.
\end{ex}

\begin{ex} \label{ex:4}
We consider \eqref{vv} only on a quadrant $U=(-\infty,0)\times (0,\infty)$.
Assume $F(p)=p^3$ for $p\in \R$, and $g(x)=\frac{2\sqrt{2}}{9} (-x)^{\nicefrac{3}{2}}$ for $x \leq 0$.
The limiting PDE is 
\begin{equation*}
\begin{cases}
u_t + (u_x)^3 = 0 \quad &\text{ in } (-\infty,0) \times (0,\infty),\\
u(0,t)=0 \quad &\text{ for } t \in (0,\infty),\\
u(x,0) = g(x) \quad &\text{ for } x \in (-\infty,0].
\end{cases}
\end{equation*} 
Then, for $0\leq t \leq 1$,
\[
u(x,t)=\left(-\frac{2x}{3}\right)^{\nicefrac{3}{2}}  (3-2t)^{-\nicefrac{1}{2}} \quad \text{ for all } x\in (-\infty,0].
\]
Numerical results are shown in Figure \ref{Fig: plots} (D). 
{Numerically, we observe that the convergence rate is $O(\ep^{3/4})$ in this example.}
\end{ex}

\begin{ex} \label{ex:5}
Assume $F(p)=\frac{1}{4}|p|^{4}$ for $p\in \R$, and $g(x)=M\min(|x|,|x-\frac{1}{2}|-\frac{1}{4})$ for $x\in \R$ and some scaling constant $M$. We  choose $M\in\{\frac{1}{4}, \frac{1}{2},  1,  2\}$ to perform our tests. Then, we can use the Hopf-Lax formula to obtain
\begin{align*}
u(x,t)= \inf_{y\in \R} \left\{g(y)+t\, L\left(\frac{x-y}{t}\right)\right\},\quad\text{where}\quad L(v):=\sup_{p\in \R} \left\{pv-F(p)\right\} = \frac{3}{4}|v|^{\nicefrac{4}{3}}
\end{align*}
for $(x,t)\in \R\times (0,\infty)$. Numerical results are shown in Figure \ref{Fig: plots} (E). 
{Numerically, we observe that the convergence rate is $O(\ep^{\nicefrac{2}{3}})$ in this example.}
\end{ex}

\subsection{A simple homogenization test}
\begin{figure}
	\begin{subfigure}{0.49\textwidth}
		\includegraphics[width=\textwidth]{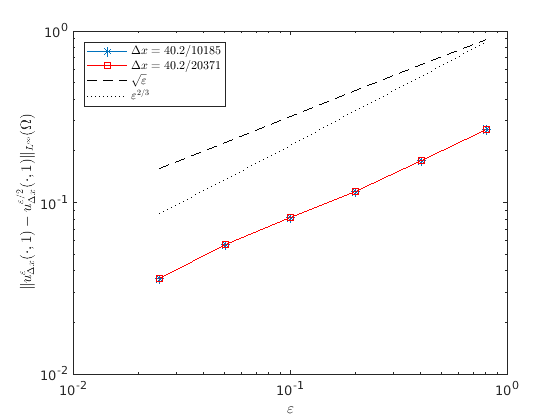}
		\subcaption{$H(y,p) \hspace{-0.05cm}=\hspace{-0.05cm} \frac{1}{2} |p|^2 \hspace{-0.05cm} +\hspace{-0.05cm} \min_{k\in \Z}\hspace{-0.05cm} |y-k|$, $y,p\in \R$.}
	\end{subfigure}
	\begin{subfigure}{0.49\textwidth}
		\includegraphics[width=\textwidth]{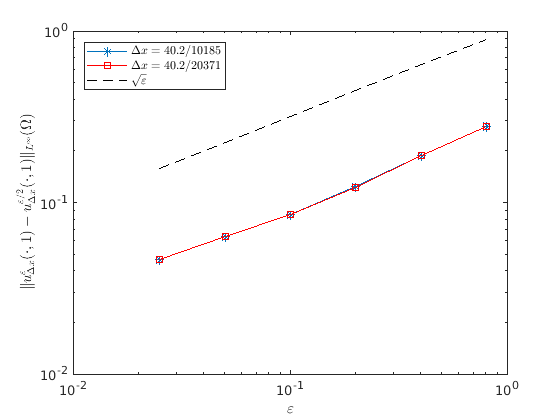}
		\subcaption{$H(y,p) \hspace{-0.05cm}=\hspace{-0.05cm} \frac{1}{4} |p|^4 \hspace{-0.05cm} +\hspace{-0.05cm} \min_{k\in \Z}\hspace{-0.05cm} |y-k|$, $y,p\in \R$.}
	\end{subfigure}
	
\begin{subfigure}{0.49\textwidth}
		\includegraphics[width=\textwidth]{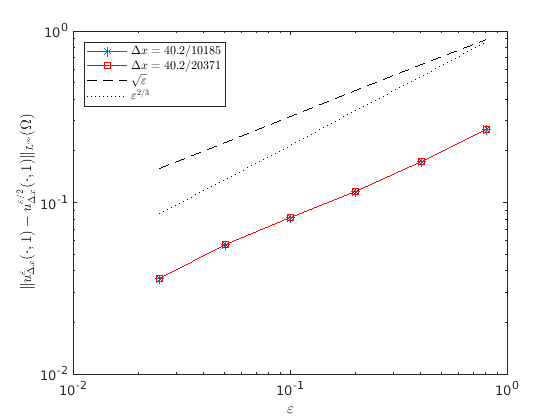}
		\subcaption{$H(y,p) \hspace{-0.05cm}=\hspace{-0.05cm} \frac{1}{2} |p|^2 \hspace{-0.05cm} +\hspace{-0.05cm} \min_{k\in \Z} \hspace{-0.05cm}|y-k|^2$\hspace{-0.05cm}, $y,p\in \R$.}
	\end{subfigure}
	\begin{subfigure}{0.49\textwidth}
		\includegraphics[width=\textwidth]{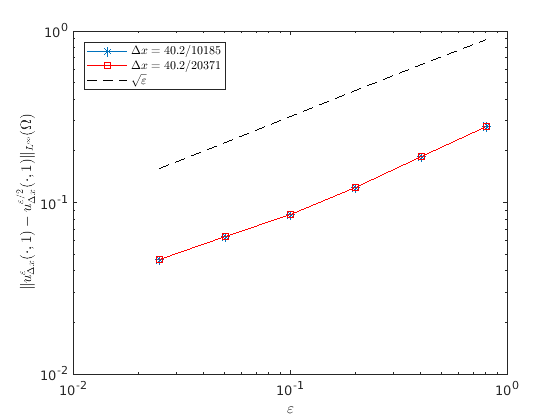}
		\subcaption{$H(y,p) \hspace{-0.05cm}=\hspace{-0.05cm} \frac{1}{4} |p|^4 \hspace{-0.05cm} +\hspace{-0.05cm} \min_{k\in \Z} \hspace{-0.05cm}|y-k|^2$\hspace{-0.05cm}, $y,p\in \R$.}
	\end{subfigure}	
	
\begin{subfigure}{0.49\textwidth}
		\includegraphics[width=\textwidth]{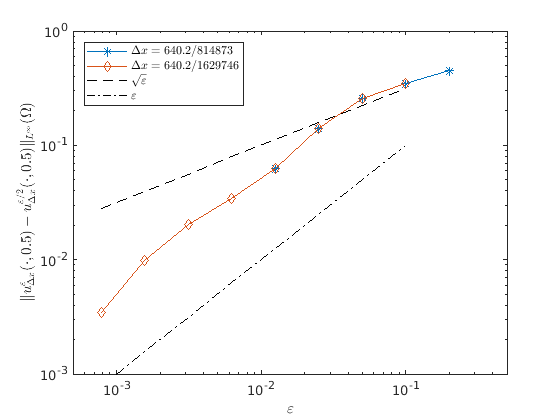}
		\subcaption{$H(y,p) = \frac{1}{2} |p|^2 + \sin(y)$ for $y,p\in \R$.}
	\end{subfigure}
	\begin{subfigure}{0.49\textwidth}
		\includegraphics[width=\textwidth]{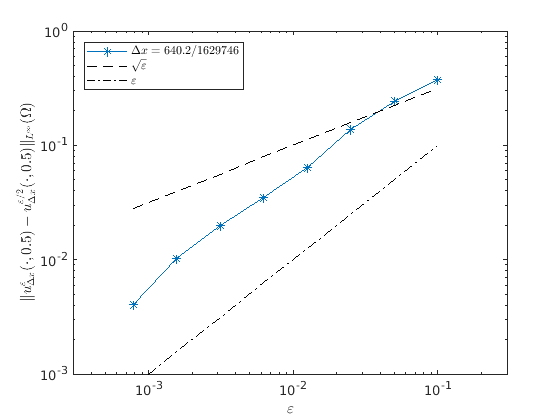}
		\subcaption{$H(y,p) = \frac{1}{4} |p|^4 + \sin(y)$ for $y,p\in \R$.}
	\end{subfigure}
	
	\caption{Illustration of $\|u^{\ep}_{\Delta x}(\cdot,T)-u^{\nicefrac{\ep}{2}}_{\Delta x}(\cdot,T)\|_{L^\infty(\Omega)}$ for Examples \ref{ex:6}--\ref{ex:11} with initial datum $g(x) = \min(|x|,|x-\frac{1}{2}|-\frac{1}{4})$ for $x\in \R$. Here, $\Omega = [-\frac{5}{2},\frac{5}{2}]$, $T = 1$ for (A)--(D), and $\Omega = [-\frac{11}{2},\frac{11}{2}]$, $T = \frac{1}{2}$ for (E)--(F).}
	\label{Fig: plotshom}
\end{figure}

Consider \eqref{eq:C-ep} in one dimension, that is,
\begin{equation*}
\begin{cases}
 u_t^\ep+H\left(\frac{x}{\ep},u^\ep_x \right)=\ep u^\ep_{xx} \qquad &\text{in} \ \R \times (0,\infty),\\
u^\ep(x,0)=g(x) \qquad &\text{on} \ \R.
\end{cases} 
\end{equation*}
We take $g(x)= \min(|x|,|x-\frac{1}{2}|-\frac{1}{4})$ for $x\in \R$, and we consider six different choices for the Hamiltonian $H$. Since the exact solution to the homogenized problem \eqref{eq:C} is unknown, we compute $\|u^{\ep}(\cdot,T)-u^{\nicefrac{\ep}{2}}(\cdot,T)\|_{L^{\infty}(\Omega)}$ for some chosen $T>0$ and computational domain $\Omega$.  
\begin{ex}\label{ex:6}
Assume $H(y,p) = \frac{1}{2}|p|^2 + \min_{k\in \Z} |y-k|$ for $y,p\in \R$. Numerical results are shown in Figure \ref{Fig: plotshom} (A). The order of convergence seems to be in $[\frac{1}{2},\frac{2}{3}]$.  
\end{ex}

\begin{ex}\label{ex:7}
Assume $H(y,p) = \frac{1}{4}|p|^4 + \min_{k\in \Z} |y-k|$ for $y,p\in \R$. Numerical results are shown in Figure \ref{Fig: plotshom} (B), and the order of convergence seems to be $\frac{1}{2}$.  
\end{ex}

\begin{ex}\label{ex:8}
Assume $H(y,p) = \frac{1}{2}|p|^2 + \min_{k\in \Z} |y-k|^2$ for $y,p\in \R$. Numerical results are shown in Figure \ref{Fig: plotshom} (C). We observe the same as for Example \ref{ex:6}.  
\end{ex}

\begin{ex}\label{ex:9}
Assume $H(y,p) = \frac{1}{4}|p|^4 + \min_{k\in \Z} |y-k|^2$ for $y,p\in \R$. Numerical results are shown in Figure \ref{Fig: plotshom} (D). We observe the same as for Example \ref{ex:7}.    
\end{ex}

\begin{ex}\label{ex:10}
Assume $H(y,p) = \frac{1}{2}|p|^2 + \sin(y)$ for $y,p\in \R$. Numerical results are shown in Figure \ref{Fig: plotshom} (E), and the order of convergence seems to be close to $1$.  
\end{ex}

\begin{ex}\label{ex:11}
Assume $H(y,p) = \frac{1}{4}|p|^4 + \sin(y)$ for $y,p\in \R$. Numerical results are shown in Figure \ref{Fig: plotshom} (F), and the order of convergence seems to be close to $1$.
\end{ex}

\section{Numerical approximation of effective Hamiltonians}\label{Sec: Num effHbar}

In this section, we would like to gain a better understanding of the effective Hamiltonian $\ol H$. Let us recall that for $p\in \R^n$, the value $\ol H(p)\in \R$ is the unique constant for which there exists a viscosity solution $v(\cdot,p)\in C(\T^n)$ to
\begin{equation*}
 \qquad H(y,p+Dv)=\ol H(p) +\Delta v \qquad \text{for} \ y\in \T^n.
\end{equation*}

\subsection{Framework}

Let us focus on a Hamilton-Jacobi-Bellman nonlinearity
\begin{align}\label{HJB Ham}
H:\T^n\times \R^n\rightarrow \R,\qquad H(y,p):= \sup_{\alpha\in\Lambda} \left\{ -b(y,\alpha)\cdot p - f(y,\alpha) \right\},
\end{align}
where $\Lambda$ is a compact metric space, $b \in C(\T^n\times \Lambda;\R^n)$, $f\in C(\T^n\times\Lambda)$, and we assume that $b = b(y,\alpha)$, $f= f(y,\alpha)$ are Lipschitz continuous in $y$, uniformly in $\alpha$.
In this setting, $H \in \Lip(\T^n \times \R^n)$ and $H = H(y,p)$ is convex in $p$.
See \cite{LXY} and the references therein for the homogenization of viscous G-equations.

\subsection{Approximation of the effective Hamiltonian}

Let $p\in \R^n$ be fixed. Our goal is to approximate the value $\ol H(p)$, and we begin by introducing approximate correctors.

\subsubsection{Approximate correctors}

For $\sigma > 0$, introducing the approximate corrector $v^{\sigma}\in C(\T^n)$ to be the unique viscosity solution to the problem
\begin{align}\label{vsig eqn}
\sigma v^{\sigma} +H(y,p+Dv^{\sigma}) =  \Delta v^{\sigma}\qquad\text{for} \ y\in \T^n,
\end{align} 
it is known that $\{-\sigma v^{\sigma}\}_{\sigma > 0}$ converges uniformly to the constant $\ol H(p)$ as $\sigma \rightarrow 0^+$; see \cite[Chapter 4]{LMT}.

\begin{lem}\label{Lmm: sig vsig}
For $\sigma>0$, let $v^\sigma$ denote the unique viscosity solution to \eqref{vsig eqn}. Then, $v^\sigma \in C^{2,\gamma}(\T^n)$ for any $\gamma\in (0,1)$.
Moreover, for any $\sigma >0$ there holds
\[
\|\sigma v^\sigma + \ol H(p) \|_{L^\infty(\T^n)} \leq C \sigma,
\]
where $C>0$ is a constant independent of $\sigma$.
\end{lem}

\begin{proof}
As $H \in \Lip(\T^n \times \R^n)$, we have that $v^\sigma\in W^{2,p}(\T^n)$ for any $p>1$; see \cite{AmCr}.
Hence, $Dv^\sigma \in C^{0,\gamma}(\T^n)$ for any $\gamma\in (0,1)$, and hence, $H(\cdot,Dv^\sigma(\cdot))\in C^{0,\gamma}(\T^n)$.
By the standard Schauder estimates, we obtain that $v^\sigma \in C^{2,\gamma}(\T^n)$.

Let $v=v(\cdot,p)\in C(\T^n)$  be a solution to the cell problem \eqref{eq:E-p}.
Then, the function $v-\|v\|_{L^\infty(\T^n)} - \frac{\ol H(p)}{\sigma}$ is a subsolution to \eqref{vsig eqn} and the function $v+\|v\|_{L^\infty(\T^n)}- \frac{\ol H(p)}{\sigma}$ is a supersolution to \eqref{vsig eqn}. By the comparison principle, we have that
\[
v-\|v\|_{L^\infty(\T^n)} - \frac{\ol H(p)}{\sigma} \leq v^\sigma \leq v+\|v\|_{L^\infty(\T^n)}- \frac{\ol H(p)}{\sigma}\quad\text{in }\T^n,
\]
and hence,
\[
\|\sigma v^\sigma + \ol H(p) \|_{L^\infty(\T^n)} \leq 2 \|v\|_{L^\infty(\T^n)}  \sigma,
\]
which completes the proof.
\end{proof}

Therefore, a natural idea is to obtain a numerical approximation of $\ol H(p)$ based on the fact that 
\begin{align*}
\ol H(p) = \lim_{\sigma\rightarrow 0^+} \int_{Y} (-\sigma v^{\sigma}), 
\end{align*}
where $Y:=(0,1)^n$, in combination with a numerical approximation $v^{\sigma}_h$ of $v^{\sigma}$ with $\|v^{\sigma}-v^{\sigma}_h\|_{L^1(Y)}\rightarrow 0$ as $h\rightarrow 0$ for $\sigma$ fixed. Let us briefly address a possible numerical approximation for $\|b\|_{\infty}$ small. To ensure strong monotonicity of the finite element schemes proposed below, we assume that 
\begin{align}\label{min disc}
\sigma > \frac{\|b\|_{\infty}^2}{4},
\end{align}
requiring a minimum discount to be imposed for the numerical scheme. {Here, we follow the idea of the small-$\delta$ method (see, e.g., \cite{Qian}) in combination with a finite element approximation of \eqref{vsig eqn}. } We note that the effective Hamiltonian can also be approximated by the large-T method; see \cite{Qian, QTY} and the references therein.
Since the large-T method and the small-$\delta$ method (see, e.g., \cite{Qian}) are mathematically equivalent, we just use the small-$\delta$ method to illustrate the new formulation for convenience.

\subsubsection{$H^1_{\mathrm{per}}(Y)$-conforming finite element approximation of \eqref{vsig eqn}} \label{Sec: H1-conf}

We have that $v^{\sigma}$ is the unique element in $H^1_{\mathrm{per}}(Y)$ such that
\begin{align*}
a(v^{\sigma},\varphi)= 0\qquad \forall \varphi\in H^1_{\mathrm{per}}(Y),
\end{align*}
where $a: H^1_{\mathrm{per}}(Y)\times H^1_{\mathrm{per}}(Y)\rightarrow \R$ is given by
\begin{align*}
a(w,\varphi):=(D w, D \varphi)_{L^2(Y)} + ( \sup_{\alpha\in\Lambda} \left\{  -b(\cdot,\alpha)\cdot Dw -g(\cdot,\alpha) \right\}, \varphi)_{L^2(Y)} + \sigma ( w, \varphi)_{L^2(Y)}
\end{align*}
with $g(\cdot,\alpha):= b(\cdot,\alpha)\cdot p + f(\cdot,\alpha)$. Indeed, assuming \eqref{min disc}, $a: H^1_{\mathrm{per}}(Y)\times H^1_{\mathrm{per}}(Y)\rightarrow \R$ is strongly monotone since for any $u_1,u_2\in H^1_{\mathrm{per}}(Y)$ and $s\in (\frac{\|b\|_{\infty}^2}{4\sigma},1)$, writing $\delta_u:= u_1-u_2$, 
\begin{align*}
a(u_1,\delta_u) - a(u_2,\delta_u) &\geq \|D \delta_u\|_{L^2(Y)}^2 + \sigma \|\delta_u\|_{L^2(Y)}^2 - ( \sup_{\alpha\in\Lambda}\ \lvert b(\cdot,\alpha)\cdot D \delta_u\rvert ,\lvert\delta_u\rvert)_{L^2(Y)}  \\&\geq (1-s)\|D \delta_u\|_{L^2(Y)}^2 + \left(\sigma - \frac{1}{4s}\|b\|_{\infty}^2\right)\|\delta_u\|_{L^2(Y)}^2 \\ &\geq C_m \|\delta_u\|_{H^1(Y)}^2.
\end{align*}
It is also quickly checked that we have the Lipschitz property 
\begin{align*}
\lvert a(u_1,\varphi)-a(u_2,\varphi)\rvert\leq C_l\|u_1-u_2\|_{H^1(Y)}\|\varphi\|_{H^1(Y)} \qquad\forall u_1,u_2,\varphi\in H^1_{\mathrm{per}}(Y).
\end{align*}
Let $V_h\subset H^1_{\mathrm{per}}(Y)$ be a closed linear subspace of $H^1_{\mathrm{per}}(Y)$. By the Browder-Minty theorem and standard conforming Galerkin arguments, there exists a unique $v^{\sigma}_h \in V_h$ such that
\begin{align}\label{Meth Hbar}
a(v^{\sigma}_h,\varphi_h) = 0\qquad \forall \varphi_h\in V_h,
\end{align}  
and we have the near-best approximation bound
\begin{align*}
\|v^{\sigma} - v^{\sigma}_h\|_{H^1(Y)}\leq \frac{C_l}{C_m}\inf_{w_h\in V_h} \|v^{\sigma} - w_h\|_{H^1(Y)}.
\end{align*}
Choosing for $V_h$ a Lagrange finite element space over a shape-regular triangulation $\mathcal{T}_h$ of $\ol Y$ with mesh-size $h>0$, consistent with the periodicity requirement, leads to a convergent method under mesh refinement. The discrete nonlinear system can be solved numerically using Howard's algorithm (see e.g., \cite{SS}). 

Introducing the approximate effective Hamiltonian
\begin{align}\label{app effH}
\ol H_{\sigma,h}(p) := \int_Y (-\sigma v^{\sigma}_h),
\end{align}
we then have that
\begin{align*}
\lvert \ol H(p) - \ol H_{\sigma,h}(p) \rvert \leq \left\lvert \ol H(p) - \int_Y (-\sigma v^{\sigma})\right\rvert + \sigma \|v^{\sigma}-v^{\sigma}_h\|_{L^1(Y)},
\end{align*}
where $\|v^{\sigma}-v^{\sigma}_h\|_{L^1(Y)}\rightarrow 0$ as $h\rightarrow 0$ and the first term on the right-hand side is of order $O(\sigma)$ by Lemma \ref{Lmm: sig vsig}.

\subsubsection{Fourth-order-type variational formulation for \eqref{vsig eqn}}

If information on second-order derivatives of $v^{\sigma}$ is desired, it is interesting to see that inspired by arguments based on Cordes-type conditions (see e.g., \cite{CSS,GSS,SS,Spr}), we can derive a fourth-order-type variational formulation for $v^{\sigma}$, allowing for the construction of $H^2$-conforming finite element schemes. Introducing $\gamma:=\frac{4\sigma}{\lvert b\rvert^2 + 4\sigma}\in C(\T^n\times \Lambda,(0,1])$, note that $v^{\sigma}$ is the Y-periodic solution to 
\begin{align*}
G[v^{\sigma}]= 0,\quad \text{where}\quad G[w]:=\sup_{\alpha\in\Lambda} \left\{ \gamma(\cdot,\alpha)\left( -\Delta w-b(\cdot,\alpha)\cdot Dw + \sigma w -g(\cdot,\alpha)\right) \right\},
\end{align*}
and $v^{\sigma}$ is the unique element in $H^2_{\mathrm{per}}(Y)$ satisfying
\begin{align*}
\tilde{a}(v^{\sigma},\varphi):= ( G[v^{\sigma}],\sigma \varphi-\Delta \varphi)_{L^2(Y)} = 0\qquad \forall \varphi\in H^2_{\mathrm{per}}(Y).
\end{align*}
Indeed, note that due to \eqref{min disc} we have that $\tilde{a}$ is strongly monotone: For any $u_1,u_2\in H^2_{\mathrm{per}}(Y)$, writing $\delta_u:=u_1-u_2$ and $\eta := \frac{4\sigma - \|b\|_{\infty}^2}{4\sigma + \|b\|_{\infty}^2} \in (0,1]$, we have 
\begin{align*}
\lvert G[u_1] &- G[u_2] - (\sigma \delta_u -\Delta \delta_u)\rvert^2 \\ &\leq \sup_{\alpha\in \Lambda} \ \lvert -(\gamma(\cdot,\alpha)-1)\Delta \delta_u - [\gamma b](\cdot,\alpha)\cdot D \delta_u + (\gamma(\cdot,\alpha)-1)\sigma \delta_u \rvert^2  
\\ &\leq (1-\eta)(\lvert \Delta \delta_u\rvert^2 + 2\sigma \lvert D \delta_u\rvert^2 + \sigma^2 \lvert \delta_u\rvert^2)
\end{align*}
almost everywhere (note $2\lvert \gamma -1\rvert^2 + \frac{1}{2\sigma}\lvert \gamma b\rvert^2 = 2-2\gamma\leq 1-\eta$), and
\begin{align*}
\|\Delta \delta_u\|_{L^2(Y)}^2 + 2\sigma \|D \delta_u\|_{L^2(Y)}^2 + \sigma^2 \|\delta_u\|_{L^2(Y)}^2 = \|\sigma\delta_u - \Delta \delta_u\|_{L^2(Y)}^2,
\end{align*}
which in combination yields
\begin{align*}
\tilde{a}(u_1,\delta_u) - \tilde{a}(u_2,\delta_u) \geq \left(1-\sqrt{1-\eta}\right) \|\sigma \delta_u -\Delta \delta_u\|_{L^2(Y)}^2.
\end{align*}
Further, $\tilde{a}$ satisfies the Lipschitz property
\begin{align*}
\lvert \tilde{a}(u_1,\varphi)-\tilde{a}(u_2,\varphi)\rvert \leq \left(\sqrt{1-\eta}+\sqrt{2}\right)\|\sigma \delta_u - \Delta \delta_u\|_{L^2(Y)}\|\sigma \varphi - \Delta \varphi\|_{L^2(Y)}.
\end{align*}
Let $V_h\subset H^2_{\mathrm{per}}(Y)$ be a closed linear subspace of $H^2_{\mathrm{per}}(Y)$. By the Browder-Minty theorem and standard conforming Galerkin arguments, there exists a unique $v^{\sigma}_h \in V_h$ such that
\begin{align*}
\tilde{a}(v^{\sigma}_h,\varphi_h) = 0\qquad \forall \varphi_h\in V_h,
\end{align*}  
and, introducing the norm $||| w ||| := \|\sigma w - \Delta w\|_{L^2(Y)}$ for $w\in H^2_{\mathrm{per}}(Y)$, we have the near-best approximation bound
\begin{align*}
|||v^{\sigma} - v^{\sigma}_h|||\leq \frac{\sqrt{1-\eta}+\sqrt{2}}{1-\sqrt{1-\eta}}\inf_{w_h\in V_h} ||| v^{\sigma} - w_h|||.
\end{align*}
Choosing for $V_h$ an Argyris or HCT finite element space over a shape-regular triangulation $\mathcal{T}_h$ of $\ol Y$ with mesh-size $h>0$, consistent with the periodicity requirement, leads to a convergent method under mesh refinement. The discrete nonlinear system can again be solved numerically using Howard's algorithm. With the observations of this subsection at hand, one can also construct mixed finite element schemes and discontinuous Galerkin finite element schemes for \eqref{vsig eqn} similarly to \cite{GSS,KS}.

\subsubsection{Numerical experiments}

For our numerical tests, we consider one linear example with known effective Hamiltonian and one nonlinear example with unknown effective Hamiltonian. For both tests, we use the method from Section \ref{Sec: H1-conf}.

\begin{figure}
	\begin{subfigure}{0.49\textwidth}
		\includegraphics[width=\textwidth]{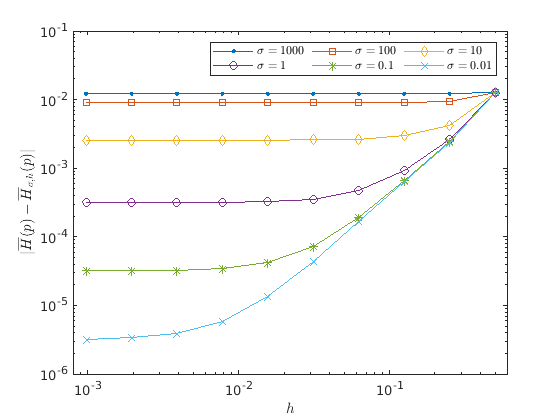}
		\subcaption{$h\mapsto\lvert \ol H(p) - \ol H_{\sigma,h}(p) \rvert$ for fixed $\sigma$}
	\end{subfigure}
	\begin{subfigure}{0.49\textwidth}
		\includegraphics[width=\textwidth]{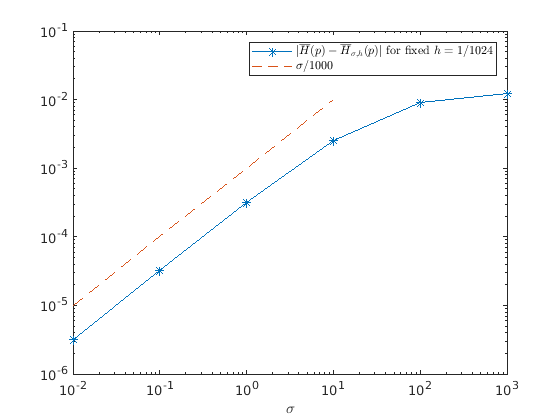}
		\subcaption{$\sigma\mapsto\lvert \ol H(p) - \ol H_{\sigma,h}(p) \rvert$ for fixed $h = 2^{-10}$}
	\end{subfigure}	
	\caption{Approximation of $\ol H(p)$ at $p = (3,1)$ for Example \ref{ex:Hbar1}.}
	\label{Fig: Hbar exp}
\end{figure}

\begin{ex}\label{ex:Hbar1}
Consider $H:\T^2\times \R^2\rightarrow \R$ given by \eqref{HJB Ham} with $n = 2$ and $\Lambda := \{0\}$. We set $b(y,\alpha):=\tilde{b}(y):=\left(\frac{1}{2\pi}\cos(2\pi y_1),0\right)$ and $f(y,\alpha):=\tilde{f}(y) =  1 + \sin(2\pi y_1)$ for $y=(y_1,y_2)\in \T^2$ and $\alpha\in \Lambda$. Our goal is to approximate the value of the effective Hamiltonian $\ol H$ at the point $p := (3,1)$, and compute the approximation error $\lvert \ol H(p) - \ol H_{\sigma,h}(p) \rvert$, where the true value can be explicitly computed as
\begin{align*}
\ol H(p) = 1+\frac{\int_0^1 \sin(2\pi t) \exp(\frac{1}{4\pi^2}\sin(2\pi t))\, \mathrm{d}t}{\int_0^1 \exp(\frac{1}{4\pi^2}\sin(2\pi s))\,\mathrm{d}s}.
\end{align*}
In our numerical experiment, we compute $\ol H_{\sigma,h}(p)$ via \eqref{Meth Hbar}--\eqref{app effH}, where we choose $V_h$ to consist of continuous $Y$-periodic piecewise affine functions on a periodic shape-regular triangulation $\mathcal{T}_h$ of $\ol Y$ into triangles with vertices $\{(ih, jh)\}_{1\leq i,j\leq N}$ where $N = \frac{1}{h}\in \N$. We choose $\sigma = 10^{-i}$ for $i\in [-3,2]\cap \Z$ and $h = 2^{-j}$ for $j\in [1,10]\cap \Z$. The results are shown in Figure \ref{Fig: Hbar exp}. Numerically, we can observe that the rate $O(\sigma)$ in Lemma \ref{Lmm: sig vsig} is optimal.
\end{ex} 

\begin{figure}
	\begin{subfigure}{0.49\textwidth}
		\includegraphics[width=\textwidth]{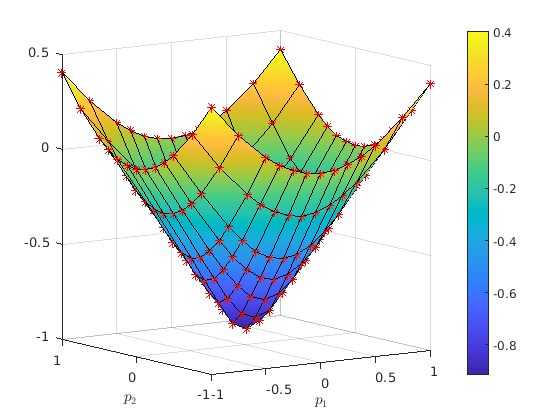}
		\subcaption{$p\mapsto \ol H_{\sigma,h}(p)$ for $\sigma = 2^{-4}$, $h = 2^{-10}$.}
	\end{subfigure}
	\begin{subfigure}{0.49\textwidth}
		\includegraphics[width=\textwidth]{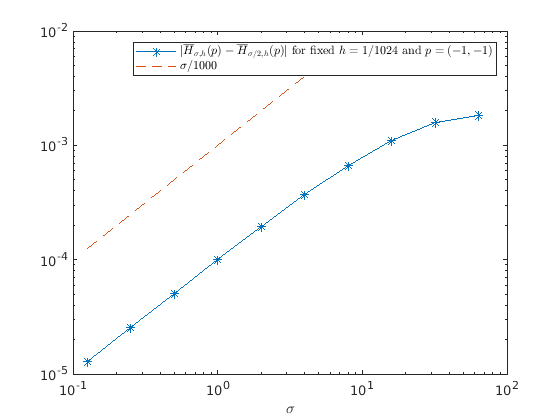}
		\subcaption{$\sigma\mapsto\lvert \ol H_{\sigma,h}(p)-\ol H_{\frac{\sigma}{2},h}(p) \rvert$ for $h = 2^{-10}$ and $p = (-1,-1)$.}
	\end{subfigure}	
	\caption{Approximation of $\ol H$ for Example \ref{ex:Hbar2}.}
	\label{Fig: Hbar exp 2}
\end{figure}

\begin{ex}\label{ex:Hbar2}
We consider $H:\T^2\times \R^2\rightarrow \R$ given by \eqref{HJB Ham} with $n = 2$ and $\Lambda := \{\alpha\in \R^2: |\alpha|\leq 1\}$. We set $b(y,\alpha):=\tilde{b}(y)+\alpha$ and $f(y,\alpha):= \tilde{f}(y)$ for $(y,\alpha)\in \T^2\times \Lambda$, where $\tilde{b}$ and $\tilde{f}$ are defined as in Example \ref{ex:Hbar1}. Note that $H(y,p) = |p|-\tilde{b}(y) \cdot p- \tilde{f}(y)$ for $(y,p)\in \T^2\times \R^2$. Our goal is to approximate the unknown effective Hamiltonian $\ol H$ on $[-1,1]^2$. To this end, we approximate $\ol H(p)$ at all points  $p$ in $S:=\{\pm 1,\pm\frac{3}{4},\pm\frac{1}{2},\pm\frac{3}{8},\pm\frac{1}{4},\pm\frac{1}{8},0\}^2$, where we chose a finer resolution around the origin. In our numerical experiment, we compute $\ol H_{\sigma,h}(p)$ via \eqref{Meth Hbar}--\eqref{app effH}, where we choose $V_h$ to consist of continuous $Y$-periodic piecewise affine functions on a periodic shape-regular triangulation $\mathcal{T}_h$ of $\ol Y$ into triangles with vertices $\{(ih, jh)\}_{1\leq i,j\leq N}$ where $N = \frac{1}{h}\in \N$. We fixed a fine mesh, i.e., $h = 2^{-10}$, and produced convergence histories with respect to $\sigma$ at each point $p\in S$. The nonlinear discrete problems were solved using Howard's algorithm. For the plot of the numerical effective Hamiltonian we used $\sigma = 2^{-4}$; see Figure \ref{Fig: Hbar exp 2} (A). An exemplary convergence history of $\lvert \ol H_{\sigma,h}(p)-\ol H_{\frac{\sigma}{2},h}(p) \rvert$ with respect to $\sigma$, for $p = (-1,-1)$, is shown in Figure \ref{Fig: Hbar exp 2} (B) and we observe the rate $O(\sigma)$, as expected.  We note that the scheme performs nicely even beyond \eqref{min disc}. 
\end{ex}

\end{document}